\documentclass[a4paper]{article}

\usepackage[utf8]{inputenc}
\usepackage{pdfpages}
\usepackage[hidelinks]{hyperref}
\usepackage{url}
\usepackage{amsthm}
\usepackage{thmtools}
\usepackage{amssymb}
\usepackage{amsmath}
\usepackage{float}
\usepackage{pdflscape}
\usepackage{natbib}
\PassOptionsToPackage{dvipsnames,x11names,table}{xcolor}

\usepackage{array}
\newcolumntype{H}{>{\setbox0=\hbox\bgroup}c<{\egroup}@{}}

\usepackage{dsfont}

\usepackage[labelformat=simple]{subcaption}

\definecolor{WtalUniGruen}{RGB}{137,186,23}

\usepackage{tikz}
\usetikzlibrary{patterns}
\usetikzlibrary{matrix}
\usetikzlibrary{graphs,quotes}
\usepackage{pgfplots}

\pgfplotsset{compat=1.17,
	every axis/.append style={
		axis lines=center,
		xlabel style={anchor=south west},
		ylabel style={anchor=south west},
		zlabel style={anchor=south west},
		tick align=outside,}
}
\usepgfplotslibrary{patchplots}
\pgfdeclareverticalshading{darkGray}{100bp}
{color(0bp)=(black!45); color(32bp)=(black!45); color(38bp)=(black!40);color(45bp)=(black!30); color(100bp)=(black!25)}

\pgfdeclareverticalshading{darkGray2}{100bp}
{color(0bp)=(black!45); color(40bp)=(black!45); color(45bp)=(black!40);color(50bp)=(black!30); color(100bp)=(black!25)}

\pgfdeclareverticalshading{midGray}{100bp}
{color(0bp)=(black!30); color(30bp)=(black!35);color(35bp)=(black!30);color(40bp)=(black!25); color(47bp)=(black!15); color(100bp)=(black!10)}

\pgfdeclareverticalshading{lightGray}{100bp}
{color(0bp)=(black!20); color(28bp)=(black!15);color(35bp)=(black!15); color(40bp)=(black!10); color(45bp)=(black!5); color(100bp)=(black!5)}

\pgfdeclareverticalshading{lightGray2}{100bp}
{color(0bp)=(black!20); color(35bp)=(black!15);color(40bp)=(black!15); color(45bp)=(black!10); color(50bp)=(black!5); color(100bp)=(white)}

\usepackage{cancel}

\usepackage[ruled,vlined,linesnumbered]{algorithm2e}
\SetKwProg{Fn}{Function}{}{}
\DontPrintSemicolon
\SetKwFor{For}{for}{}{}
\SetKwFor{While}{while}{}{}

\usepackage{rotating} \usepackage{booktabs} \usepackage{mathtools}

\newtheorem{defi}{Definition}
\newtheorem{theorem}[defi]{Theorem}
\newtheorem{lem}[defi]{Lemma}

\newtheorem{rem}[defi]{Remark}

\newtheorem{example}[defi]{Example}

\newcommand {\Z}{\mathds{Z}}

\newcommand {\R}{\mathds{R}}

\DeclareMathOperator{\interior}{int} \DeclareMathOperator{\clos}{cl}
\DeclareMathOperator{\rank}{rank}
\DeclareMathOperator{\Nset}{N}

\DeclareMathOperator{\lex}{lex}

\newcommand{\dS}{S} \newcommand{\zulM}{X} \newcommand{\xone}{x'} \newcommand{\xtwo}{\hat{x}}  \newcommand{\yone}{y'} \newcommand{\ytwo}{\hat{y}}  

\newcommand{\ones}{\boldsymbol{1}}
\newcommand{\zeros}{\boldsymbol{0}}
\newcommand{\w}{\omega}
\newcommand{\g}{\gamma}
\newcommand{\og}{(\omega,\gamma)}
\newcommand{\pareto}{P} \newcommand{\vR}{\nu} \newcommand{\VR}{\mathcal{V}^{\og}} \newcommand{\geqqV}{\geqq}\newcommand{\gV}{>}\DeclareMathOperator{\preceqqnu}{\preceqq\!} \newcommand{\preceqnu}{\preceq} \newcommand{\tailC}{C_{{\preceqqsmall}_{\og}}}

\DeclareMathOperator{\hcone}{hcone} \DeclareMathOperator{\vcone}{vcone} 

\newcommand{\preceqq}{\;\raisebox{-0.12em}{\scalebox{1.15}{\ensuremath{\mathrel{\substack{\prec\\[-.15em]=}}}}}\;}
\newcommand{\preceqqo}{\;\raisebox{-0.12em}{\scalebox{1.15}{\ensuremath{\mathrel{\substack{\prec\\[-.15em]=}_o}}}}\;}
\newcommand{\preceqqsmall}{\;\raisebox{-0.12em}{\scalebox{1.0}{\ensuremath{\mathrel{\substack{\prec\\[-.15em]=}}}}}\;}

\makeatletter
\addtotheorempostfoothook{\global\everypar{{\setbox\z@\lastbox}\global\everypar{}}}
\makeatother

\defcitealias{geovelo}{\slshape geovelo}

\usepackage{authblk}

\begin{document}

\title{Adapting Polyhedral Dominance Cones to Ordinal Preference Structures}
\author{Kathrin Klamroth}
\author{Michael Stiglmayr}
\author{Julia Sudhoff Santos}
\affil{School of Mathematics and Natural Sciences, \\University of Wuppertal, Wuppertal, Germany}
\date{}
\maketitle

\begin{abstract}
	In combinatorial optimization, ordinal costs can be used to model the quality of elements whenever numerical values are not available. When considering, for example, routing problems for cyclists, the safety of a street can be ranked in ordered categories like safe (separate bike lane), medium safe (street with a bike lane) and unsafe (street without a bike lane). However, ordinal optimization may suggest unrealistic solutions with huge detours to avoid unsafe street segments.
    
    In this paper, we investigate how partial preference information regarding the relative quality of the ordinal categories can be used to improve the relevance of the computed solutions. By introducing preference weights which describe how much better a category is at least \emph{or} at most, compared to the subsequent category, we enlarge the ordinal dominance cone. This leads to a smaller set of alternatives, i.\,e., of ordinally efficient solutions.
    The corresponding weighted ordinal ordering cone is a polyhedral cone and we provide descriptions via its extreme rays and via its facets. The latter implies a linear transformation to an associated multi-objective optimization problem w.\,r.\,t.\ Pareto dominance. This paves the way for the application of standard multi-objective solution approaches. We demonstrate the usefulness of the weighted ordinal ordering cone by investigating a safest path problem with different preference weights. Moreover, we investigate the interrelation between the weighted ordering cone to standard dominance concepts of multi-objective optimization, like, e.\,g., Pareto dominance, lexicographic dominance and weighted sum dominance. 

    \medskip\noindent
    Keywords: multiple objective programming, polyhedral ordering cones, ordinal objective functions, combinatorial optimization
\end{abstract}

\section{Introduction}\label{sec:intro}
In many applications, particularly in the context of combinatorial optimization problems (like knapsack problems, shortest path problems, and many others), a precise quantification of costs may be difficult or even impossible. Examples for criteria that are usually hard to assess are the safety of paths (e.\,g., for cyclists), the quality of clothes, sustainability, the healthiness of food, or the service of a hotel. Hotels, for example, may be classified by a system of up to five stars. However, these stars generally cannot be translated into numerical values. Indeed, we cannot say whether two stars are twice as good as one star; all we know is that two stars are (usually) better than one star. The categories are also not additive, i.\,e., spending one night in a hotel with one star and the second night in a hotel with three stars is, in general, not equivalent to spending two nights in a two star hotel. Nevertheless, the star classification of hotels implies a ranking: Hotels with more stars are preferred over hotels with fewer stars. 

Cost structures that are based on ordered categories rather than precise numerical values are called \emph{ordinal}. We consider ordinal shortest path problems throughout this paper to illustrate the concepts. Assuming that edges in a path are ranked in two categories: safe (green), and unsafe (red). Then two solutions can be compared by counting the number of elements in the respective categories. Indeed, a rational decision maker prefers one green and one red (i.\,e., $(g,r)$) over two red (i.\,e., $(r,r)$), and also over two green and one red (i.\,e., $(g,g,r)$). However, whether two green, $(g,g)$, are preferred over one red, $(r)$, depends on the preferences of the decision maker. 
See \citet{koeksalan88} for an early reference, \citet{SCHAFER20,SCHAFER:knapsack}  for applications to shortest paths and knapsack problems, and \citet{Klamroth2023Ordinal} for a theoretical analysis of ordinal cones.

\begin{figure}[H]
		\centering
		\subcaptionbox{One red versus nine green \label{fig:bspIntro1}}
						[.47\linewidth]{

\begin{tikzpicture}[scale=1,every edge quotes/.append style={font=\footnotesize}]\node[draw,circle,minimum size=0.5cm] (1) at (0,0)[] {$\scriptstyle t$};
\node[draw,circle,minimum size=0.5cm] (2) at (1.1,0) {};
\node[draw,circle,minimum size=0.5cm] (3) at (2.2,0) {};
\node[draw,circle,minimum size=0.5cm] (4) at (3.3,0) {};
\node[draw,circle,minimum size=0.5cm] (5) at (4.4,0) {};
\node[draw,circle,minimum size=0.5cm] (6) at (0,1)[] {$\scriptstyle s$};
\node[draw,circle,minimum size=0.5cm] (7) at (1.1,1) {};
\node[draw,circle,minimum size=0.5cm] (8) at (2.2,1) {};
\node[draw,circle,minimum size=0.5cm] (9) at (3.3,1) {};
\node[draw,circle,minimum size=0.5cm] (10) at (4.4,1) {};

\graph {
    (2) ->[-latex,densely dotted,thick,swap,draw=black!50!green] (1);
    (3) ->[-latex,densely dotted,thick,swap,draw=black!50!green] (2);
    (4) ->[-latex,densely dotted,thick,swap,draw=black!50!green] (3);
    (5) ->[-latex,densely dotted,thick,swap,draw=black!50!green] (4);
    (10) ->[-latex,densely dotted,thick,swap,draw=black!50!green] (5);
    (6) ->[-latex,densely dotted,thick,swap,draw=black!50!green] (7);
    (7) ->[-latex,densely dotted,thick,swap,draw=black!50!green] (8);
    (8) ->[-latex,densely dotted,thick,swap,draw=black!50!green] (9);
    (9) ->[-latex,densely dotted,thick,swap,draw=black!50!green] (10);

    (6) ->[-latex,thick,swap,draw=black!40!red] (1);
};
\end{tikzpicture}

 }
            \subcaptionbox{Six green versus four red \label{fig:bspIntro2}}
						[.47\linewidth]{

\begin{tikzpicture}[scale=1,every edge quotes/.append style={font=\footnotesize}]\node[draw,circle,minimum size=0.5cm] (1) at (0,0)[] {};
\node[draw,circle,minimum size=0.5cm] (2) at (1.1,0) {};
\node[draw,circle,minimum size=0.5cm] (3) at (2.2,0) {};
\node[draw,circle,minimum size=0.5cm] (4) at (3.3,0) {};
\node[draw,circle,minimum size=0.5cm] (5) at (4.4,0) {};
\node[draw,circle,minimum size=0.5cm] (6) at (0,1)[] {$\scriptstyle s $};
\node[draw,circle,minimum size=0.5cm] (7) at (1.1,1) {};
\node[draw,circle,minimum size=0.5cm] (8) at (2.2,1) {};
\node[draw,circle,minimum size=0.5cm] (9) at (3.3,1) {};
\node[draw,circle,minimum size=0.5cm] (10) at (4.4,1) {$\scriptstyle t$};

\graph {
    (1) ->[-latex,densely dotted,thick,swap,draw=black!50!green] (2);
    (2) ->[-latex,densely dotted,thick,swap,draw=black!50!green] (3);
    (3) ->[-latex,densely dotted,thick,swap,draw=black!50!green] (4);
    (4) ->[-latex,densely dotted,thick,swap,draw=black!50!green] (5);
    (6) ->[-latex,thick,swap,draw=black!40!red] (7);
    (7) ->[-latex,thick,swap,draw=black!40!red] (8);
    (8) ->[-latex,thick,swap,draw=black!40!red] (9);
    (9) ->[-latex,thick,swap,draw=black!40!red] (10);

    (6) ->[-latex,densely dotted,thick,swap,draw=black!50!green] (1);
    (5) ->[-latex,densely dotted,thick,swap,draw=black!50!green] (10);
};
\end{tikzpicture}

 }
		\caption{Illustration of unfavorable trade-offs between $s$-$t$-paths that may be used to extract additional preference information, e.\,g., ``one red is preferred over nine green'' (left) or ``six green are preferred over four red'' (right).}
		\label{fig:BspIntro}
\end{figure}
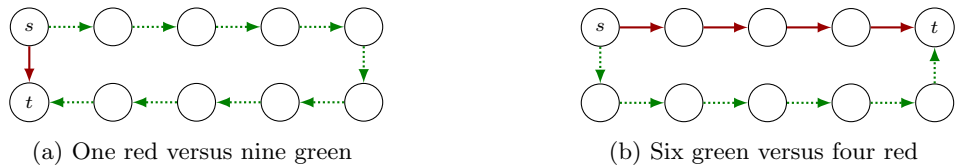

In practice, the ordinal dominance concept may lead to unfavorable solutions. Indeed, while in the case of  $(g,g)$ and $(r)$  both paths have their practical relevance, this is less clear when comparing a path with a large number of green edges with a path with only one red edge (e.\,g., $(g,g,g,g,g,g,g,g,g)$ versus $(r)$, see Figure~\ref{fig:bspIntro1}). Indeed, many (if not most) decision makers would take the risk of one red edge to avoid an overly long detour. In such cases,  additional preference information of the form ``one red is preferred over nine green'' may be available in practice. Similarly, when only a small detour is needed in order to avoid a relatively high risk (e.\,g., $(g,g,g,g,g,g)$ versus $(r,r,r,r)$, see Figure~\ref{fig:bspIntro2}), additional preference information of the form ``six green are preferred over four red'' may be available. 
We emphasize that such preference information does generally not induce a complete quantification of ordinal costs. It does, however, help to avoid somewhat extreme and unfavorable solutions.
As a consequence, the additional preference information leads to fewer non-dominated solutions. It hence bridges between purely ordinal costs and classical quantitative costs. 

In this paper, we analyze how additional preference information changes the ordinal ordering cone and discuss the consequences for appropriate solution strategies. The analysis leads to a general discussion of polyhedral ordering cones and their interrelation with standard dominance concepts like, e.\,g., Pareto dominance, lexicographic dominance, and dominance w.\,r.\,t.\ specific weighted sum scalarizations. We mainly use path planning problems with two or three safety categories to illustrate the concepts. Our results emphasize the power of preference information, and they highlight its impact on the associated nondominated sets and, indirectly, also on the required computational time to find all reasonable solution alternatives.

Our work is related to prior work on ordinal costs as well as to the topic of ordering cones in multi-objective optimization. Ordinal costs have been considered in the context of many different applications, and they have attracted increased attention recently. Early references are, for example, \cite{Bartee1971Problem} and \citet{koeksalan88}. Most articles focus on combinatorial optimization problems like, e.\,g., the shortest path problem with ordinal costs and partially ordered edges in \cite{SCHAFER20} and \cite{Bossong1999Minimal}, respectively, distributing indivisible items among agents in \cite{Fishburn1996Binary,Brams2003Fair,Brams2005Efficient,Bouveret2010Fair}, the ordinal knapsack problem in \cite{Delort2011Committee,SCHAFER:knapsack}, the minimum spanning tree problem on graphs with partially ordered edges in \cite{Schweigert1999Ordered}, and ordinal costs in matroid optimization, see \cite{Klamroth2023Multi}. In relation to multi-objective optimization, ordered categories were also used to reduce the set of Pareto-non-dominated points, see \cite{OMahony2013Sorted}.

Different notions of ordinal dominance are suggested and employed in the literature on ordinal costs. Common concepts are, for example, concepts based on sorted category vectors, injective mappings, numerical representations, or a multi-objective reformulation. Many of these concepts turn out to be equivalent and allow for a representation w.\,r.\,t.\ \emph{ordinal dominance cones}. Ordinal dominance cones are polyhedral cones that contain the positive orthant. We refer to \cite{Klamroth2023Ordinal}  for a thorough analysis and rigorous equivalence proofs in the context of  combinatorial optimization problems. 

The above discussion indicates that there is a close interrelation between ordinal preferences and cone-based preferences, which are more common in the field of
multi-objective optimization. Since the literature in this field is immense, we focus on selected results that are particularly relevant for ordinal optimization. For a general introduction to multi-objective optimization, see the books of \cite{Ehrg05} and \cite{miettinen1998}. 
An in-depth treatment of cones, including theoretical analyses, can be found, e.\,g., in \citet{goepfert23,khan2015,engau2007domination,Gerth90,yu74cone}. 

The incorporation of preference information into a cone-based dominance structure is reviewed, e.\,g., in \citet{karsu2013} and \citet{nasrabadietalsurvey}. While the former focuses more on the usage of convex cones in the context of interactive methods, the latter provides an exhaustive survey and categorization of  different approaches.
\citet{huntwiecek2003} and \cite{wiecek2007} discuss  the interrelation between the preferences of a decision-maker and dominance cones in general. \cite{Hunt2010Relative} formulate pairwise preference relations using allowable and directional trade-offs to model the preference structure between objective functions. This lead to adapted polyhedral ordering cones for which both representations are determined: their extreme rays and their describing inequalities. \citet{Kaddani2016Partial} and \cite{kaddani:tel-01827021} suggest translated ordering cones to represent decision-makers preferences, while  \cite{Ramesh1989Preference} focuses particularly on the incorporation of preferences  in the context of multi-objective integer programming problems.
In \cite{Ahn2014Establishing}, a polyhedral dominance cone is introduced based on intervals of weighting vectors that model different interrelations of the objective functions. This complementary approach will be reviewed in detail and used for the analysis of the weighted ordinal cone in Section~\ref{subsec:descr}.

The remainder of this paper is organized as follows. In Section~\ref{sec:basics} we introduce an optimality concept for ordinal costs, extend this by weights and define the weighted ordinal optimization problem. In Section~\ref{sec:cones}, we recall some basics on binary relations, ordering cones, and cone dominance. This is used in Section~\ref{sec:weightedCone} to characterize the weighted ordinal dominance cone. Moreover, the description of the cone through its facets is given. This is used to transform the problem into a corresponding multi-objective optimization problem. In Section~\ref{sec:numRes} we use this transformation to compute some ordinal shortest paths with respect to different weights. We conclude this paper with a summary in Section~\ref{sec:conclusion}.

\section{Optimality Concepts for Ordinal Costs}\label{sec:basics}

We consider combinatorial optimization problems with an ordinal counting objective function \eqref{eq:OCOP}, which is defined as
\begin{equation}\label{eq:OCOP}\tag{OCOP}
	\begin{array}{rl}
		\min_{\preceqnu_{o}} & c(x)\\
		\text{s.\,t.} & x\in \zulM.
	\end{array}
\end{equation}
Thereby, the set of feasible solutions $\zulM$ is a subset of the power set of a finite set $\dS$ of elements, i.\,e., $\zulM\subseteq2^\dS$. The standard ordering relation $\preceqnu_{o}$ for ordinal costs is defined in Definition~\ref{def:ordOpt} below.  A set of $K$ ordered categories $\mathcal{C}=\{\eta_1,\ldots,\eta_K\}$ is given and every element of $\dS$ is assigned to a category by a mapping $o:\dS\to \mathcal{C}$. The counting vector $c:\zulM \to\Z_\geqq^K$ 
equals in its \(i\)-th component to the number of elements in $x$ which are in category $\eta_i$, i.\,e., 
$c_i(x)=\vert\{e\in x \colon o(e)=\eta_i\}\vert$. In some situations, it is not meaningful to just count the number of elements in each category. In path problems, for example, each edge is associated with a positive real value, namely, the length of the edge. Let $l:\dS\to\R_>$ be a function that assigns a positive real value to each element, then the counting vector is defined as $c_i(x)=\sum_{\{e\in x \colon o(e)=\eta_i\}} l(e)$.

In ordinal optimization it is often assumed that a category $\eta_i$ is strictly preferred over category $\eta_{i+1}$, denoted as $\eta_i\prec_o\eta_{i+1}$  for all $i=1,\dots,K-1$, see, e.\,g., \cite{SCHAFER:knapsack,Klamroth2023Multi}. Ordinal optimality can then be defined based on the notion of \emph{strict standard numerical representations}, i.\,e., vectors $\nu\in \mathcal{V}_o \coloneqq \{\nu\in\R^K_> \colon \nu_{i}<\nu_{i+1}$ for all $i=1,\dots,K-1\}$, where $\mathcal{V}_o$ denotes the set of all strict standard numerical representations:

 \begin{defi}[cf.~\citealp{SCHAFER:knapsack,Klamroth2023Ordinal}]\label{def:ordOpt}
Let $\xone,\xtwo\in\zulM$ be two feasible solutions. Then,
\begin{enumerate}
	\item $\xone$ \emph{weakly $(\ones,\zeros)$-ordinally dominates} $\xtwo$ and $c(\xone)$ \emph{weakly $(\ones,\zeros)$-ordinally dominates} $c(\xtwo)$, denoted by $\xone\preceqqo \xtwo$ and $c(\xone)\preceqqo c(\xtwo)$, respectively, if and only if for \emph{every} $\nu\in\mathcal{V}_o$, it holds that $\nu^\top c(\xone)\leq \nu^\top c(\xtwo)$.
	\item $\xone$ \emph{$(\ones,\zeros)$-ordinally dominates} $\xtwo$ and $c(\xone)$ \emph{$(\ones,\zeros)$-ordi\-nally dominates} $c(\xtwo)$, denoted by $\xone\preceqnu_{o} \xtwo$ and $c(\xone)\preceqnu_{o} c(\xtwo)$, respectively, if and only if  $\xone$ weakly $(\ones,\zeros)$-ordinally dominates $\xtwo$ and there exists $\nu^*\in\mathcal{V}_o$ such that $\nu^{*\top} c(\xone)< \vR^{*\top} c(\xtwo)$.
	\item $x^*\in\zulM$ is called \emph{$(\ones,\zeros)$-ordinally efficient}, if there does not exist an $x\in\zulM$ such that $x\preceqnu_{o} x^*$.
	\item $c(x^*)$ is called \emph{$(\ones,\zeros)$-ordinally non-dominated outcome vector} of problem~\eqref{eq:OCOP}, if $x^*$ is $(\ones,\zeros)$-ordinally efficient.
\end{enumerate}
\end{defi}
The notation of $(\ones,\zeros)$-ordinal dominance is used to later distinguish between standard ordinal dominance and weighted ordinal dominance, see Definition~\ref{def:ordOpt-wg} below. See also Remark~\ref{rem:merge} below for arguments showing that this notation is consistent.

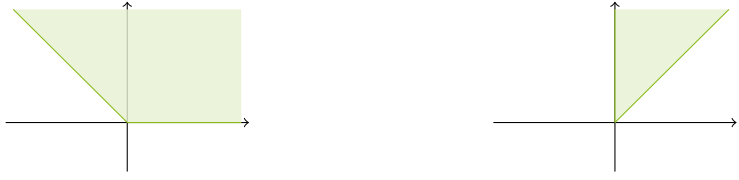
\begin{figure}
	\centering
        \subcaptionbox{$C_o\coloneqq (C_o^*)^*$: Standard ordinal cone \label{fig:ordCone2}}
						[.4\linewidth]{\begin{tikzpicture}[scale=0.65]
\draw[->] (-3.2,0) -- (3.2,0) node[right] {};
	\draw[->] (0,-1.3) -- (0,3.2) node[above] {};
 
\fill[fill=WtalUniGruen!20,fill opacity=0.8] (3,0) -- (3,3) -- (-3,3) -- (0,0)  -- cycle;
        \draw[WtalUniGruen] (0,0) -- (-3,3) node[right] {};
        \draw[WtalUniGruen] (0,0) -- (3,0) node[right] {};

\end{tikzpicture} }
        \hspace{0.7cm}
        \subcaptionbox{$C_o^*=\clos(\mathcal{V}_o)$: Dual standard ordinal cone \label{fig:DualOrdCone2}}
						[.4\linewidth]{\begin{tikzpicture}[scale=0.65]
\draw[->] (-3.2,0) -- (3.2,0) node[right] {};
	\draw[->] (0,-1.3) -- (0,3.2) node[above] {};
 
\fill[fill=WtalUniGruen!20,fill opacity=0.8] (3,3) -- (0,3) -- (0,0)  -- cycle;
        \draw[WtalUniGruen] (0,0) -- (3,3) node[right] {};
        \draw[WtalUniGruen] (0,0) -- (0,3) node[right] {};

\end{tikzpicture} }
	\caption{Standard ordinal cone $C_o$ (left) and the cone $\clos(\mathcal{V}_o)=C_o^*$ induced by the set of strict standard numerical representations  (right), which is the dual of the standard ordinal cone $C_o$.
}
		\label{fig:ordCone}
\end{figure}

Note that the closure of the set of strict standard numerical representations $\mathcal{V}_o$ defines a convex cone in $\R^K$ as follows: $C^*_o\coloneqq \clos(\mathcal{V}_o) = \{\nu\in\R^K_{\geqq} \colon 0\leq\nu_1\leq\cdots\leq\nu_K\}$. Indeed, if $\nu\in C^*_o$, then $\lambda\,\nu\in C^*_o$ also holds for all $\lambda>0$. We will see later in Section~\ref{sec:cones} that the cone $C^*_o$ is actually the dual cone of the ordering cone $C_o$ that represents standard ordinal dominance relation in $\R^K$. The two cones are illustrated in Figure~\ref{fig:ordCone}.

In the following, we generalize the concept of standard ordinal dominance by specifying the interrelation of consecutive categories by including additional preference information through marginal weights.
Hence, we assume that two weight vectors $\w,\gamma\in\R^{K-1}$ are given, with $\w_i\geq1$ and $\g_i\geq0$ for $i=1,\dots,K-1$, which model the additional preference information that interrelates consecutive categories. The weights $\w_i$ define the number of elements of a category $\eta_i$ which are (strictly) preferred over \emph{one} element of category $\eta_{i+1}$ for $i=1,\dots,K-1$, denoted as $\w_i\eta_i\preceqnu_{\og}\eta_{i+1}$ ($\w_i\eta_i\prec_{\og}\eta_{i+1}$, respectively). For $\g_i>0$, the fraction $1/\g_i$ represents the number of elements of a category $\eta_i$ which are (strictly) worse than \emph{one} element of category $\eta_{i+1}$ for $i=1,\dots,K-1$, denoted as $\g_i\eta_{i+1}\preceqnu_{\og}\eta_{i}$ ($\g_i\eta_{i+1}\prec_{\og}\eta_{i}$, respectively). If $\g_i=0$, then no set of elements of category $\eta_i$ is dominated by an element of category $\eta_{i+1}$, i.\,e., no additional preference information is provided. Moreover, if $\g_i=0$ and $\w_i=1$, we obtain the classical ordinal preference structure.
We obtain the following weighted ordinal optimization problem
\begin{equation}\label{eq:WOOP}\tag{$\text{WOOP}_{\og}$}
	\begin{array}{rl}
		\min_{\preceqnu_{\og}} & c(x)\\
		\text{s.\,t.} & x\in \zulM,
	\end{array}
\end{equation}
based on the weighted ordinal dominance relation $\preceqnu_{\og}$, which is formally defined in Definition~\ref{def:ordOpt-wg} below. First, we consider the following example to illustrate the meaning of the weights.

\begin{example}\label{ex:2}
    We consider a weighted ordinal optimization problem \ref{eq:WOOP}  with three categories $\eta_1,\eta_2,\eta_3$. For $\w_1=\w_2=1$ and $\g_1=\g_2=0$ the category $\eta_i$ is strictly preferred over category $\eta_{i+1}$ for $i=1,2$. This corresponds to the standard ordinal dominance relation as, for example, considered in \cite{SCHAFER:knapsack}, \cite{Klamroth2023Multi}, \cite{Klamroth2023Ordinal}. In this case, a set of at least two elements in the best category is incomparable to an element of the second category. 
    
    If we know, e.\,g., that ten elements of the best category are still better than one element of the second category, but we can not compare eleven elements of the best category with one element of the second category, then we have to choose $\w_1=10$. On the other hand, if we know that ten elements of the best category are worse than one element of the second category, but we can not compare nine elements of the best category with one element of the second category, then we have to choose $\g_1=\frac{1}{10}$. If we set  $\w_1=10$ and $\g_1=\frac{1}{10}$ at the same time, then ten elements of the best category are exactly as good as one element of the second category. In this case, we have no strict preference structure, see Definition~\ref{def:generalNumRep} below.

    If we have that $w_1=2$ elements of category $\eta_1$ are (strictly) preferred over one element of category $\eta_2$ and $w_2=3$ elements of category $\eta_2$ are (strictly) preferred over one element of category $\eta_3$ then it follows that $w_1\,w_2=6$ elements of category $\eta_1$ are (strictly) preferred over one element of category $\eta_3$, see Definition~\ref{def:generalNumRep} below.
\end{example}

Example~\ref{ex:2} suggests that meaningful weight vectors in the context of ordinal optimization should satisfy $\w \geqq 1$. Since we will also consider weighted preferences in a more general context, we define a slightly more general set of feasible weight vectors $\bar{\Omega}$ already at this point:
\begin{equation*}
    \bar{\Omega}\coloneqq\bigl\{(\w,\g)\in\R^{K-1}\times \R^{K-1}: 0\leqq\w,\ 0\leqq\g \text{ and } \w_i\,\g_i\leq 1 \;\forall i=1,\dots,K-1\bigr\}
\end{equation*}
Nevertheless, in the context of ordinal optimization, the reader can safely assume that $\w\geqq 1$.

Based on these considerations, we formally introduce weighted numerical representations and strict weighted numerical representations to clarify the difference.
\begin{defi}\label{def:generalNumRep}
    A vector $\vR\in\R^K_>$ is called a \emph{strict $\og$-numerical representation} with respect to given weights $\og\in\bar{\Omega}$ if
    \begin{align} 
	 & \w_i\, \nu_{i}<\nu_{i+1}\text{ for all }i\in\{1,\dots,K-1\}\quad \text{and} \label{eq:1-strict}\\ 
    & \nu_{i}>\g_i\, \nu_{i+1}\text{ for all }i\in\{1,\dots,K-1\}.  \label{eq:2-strict}
    \end{align}
Similarly, a vector $\vR\in\R^K_\geqq$ is called a \emph{$\og$-numerical representation} with respect to given weights $\w,\g\in\R_\geqq^{K-1}$ if 
    \begin{align}
	 & \w_i\, \nu_{i}\leq \nu_{i+1}\text{ for all }i\in\{1,\dots,K-1\}\quad \text{and} \label{eq:3}\\
    &\nu_{i}\geq\g_i\,\nu_{i+1}\text{ for all }i\in\{1,\dots,K-1\}. \label{eq:4}
\end{align}
    We denote by $\VR_\gV\subseteq\R^K_>$ the \emph{set of all strict $\og$-numerical representations} and with $\VR_\geqqV\subseteq \R^K_\geqq$ the \emph{set of all $\og$-numerical representations} for a given number of categories $K$ and given weights $\w,\g\in\R_\geqq^{K-1}$.
\end{defi}

To simplify the notation in the context of Definition~\ref{def:generalNumRep} throughout this paper, we define the vectors 
\begin{equation}\label{eq:u}
u^i\coloneqq (0,\dots,0,-\w_i,1,0,\dots,0)^\top\in\R^K
\end{equation}
with $-\w_i$ in its $i$-th component and 
\begin{equation}\label{eq:g} 
g^i\coloneqq (0,\dots,0,1,-\g_i,0,\dots,0)^\top\in\R^K
\end{equation}
with $-\g_i$ in its $(i+1)$-st component for all $i=1,\dots,K-1$. 
Then, \(\nu^\top u^i \geq 0\) and \(\nu^\top g^i \geq 0\) for all \(i=1,\ldots, K-1\) if and only if \(\nu\in\R^K_{\geqq}\) satisfies equations \eqref{eq:3} and \eqref{eq:4}, i.\,e., if and only if \(\nu\in\VR_{\geqq}\).

\begin{rem}\label{rem:merge}
    Note that it holds $\VR_\gV\subset \VR_\geqqV$. Moreover, the standard ordinal dominance structure is obtained for $\w=\ones\coloneqq(1,\dots,1)^\top\in\R^{K-1}$ and $\g=\zeros\coloneqq(0,\dots,0)^\top\in\R^{K-1}$, i.\,e., $\mathcal{V}_o=\mathcal{V}^{(\ones,\zeros)}_>$.
Obviously, it has to hold that $\w_i\,\g_i<1$ to distinguish the categories $\eta_i$ and $\eta_{i+1}$. Otherwise, if $\w_i\,\g_i=1$ holds then $w_i$ elements of category $\eta_i$ are exactly as good as one element of category $\eta_{i+1}$ and we have to use $\og$-numerical representations such that it holds 
$$\w_i\, \nu_{i}= \nu_{i+1} =\frac{1}{\g_i}\,\nu_{i}.$$ 
Consequently, in this case the weighted ordinal optimization problem~\ref{eq:WOOP} can be reformulated as a problem with $K-1$ categories by merging the categories \(\eta_i\) and \(\eta_{i+1}\). Thus, we assume in the following that $\w_i\,\g_i<1$ for all $i=1,\dots,K-1$ and denote this by the set $$\Omega\coloneqq\{\og\in\bar{\Omega}: \w_i\,\g_i<1 \text{ for all }i=1,\dots,K-1\}.$$
\end{rem}

If $\w=\g=\zeros$ the categories are independent, and hence an arbitrary positive natural number can be assigned to each category.
 
The inequalities~\eqref{eq:1-strict} and \eqref{eq:2-strict} as well as inequalities~\eqref{eq:3} and \eqref{eq:4} imply for $i'<\hat{\iota}$ that
\begin{align*}
	&\biggl(\prod_{j=i'}^{\hat{\iota}-1}\w_j\biggr)\,\nu_{i'}<\nu_{\hat{\iota}},\quad  \biggl(\prod_{j=i'}^{\hat{\iota}-1}\frac{1}{\g_j}\biggr)\,\nu_{i'}>\nu_{\hat{\iota}}\quad\text{and}\\
 &\biggl(\prod_{j=i'}^{\hat{\iota}-1}\w_j\biggr)\,\nu_{i'}\leq \nu_{\hat{\iota}},\quad  \biggl(\prod_{j=i'}^{\hat{\iota}-1}\frac{1}{\g_j}\biggr)\,\nu_{i'}\geq \nu_{\hat{\iota}},
	\end{align*}
 respectively.

For a given $\og$-numerical representation $\vR\in\VR_\geqqV$, we define the \emph{numerical value of a feasible solution} $x=\{e_1,\dots,e_n\}\in\zulM$ w.r.t. $\vR$ \citep[cf.][]{SCHAFER:knapsack,Klamroth2023Ordinal} as 
$\vR^\top\,c(x)$.
This concept is used to define the dominance structure for problem~\eqref{eq:WOOP} as a generalization of Definition~\ref{def:ordOpt}.

\begin{defi}\label{def:ordOpt-wg}
Let $\xone,\xtwo\in\zulM$ be feasible solutions and let $\og\in\bar{\Omega}$. Then,
\begin{enumerate}
	\item $\xone$ \emph{weakly $\og$-ordinally dominates} $\xtwo$ and $c(\xone)$ \emph{weakly $\og$-ordinally dominates} $c(\xtwo)$, denoted by $\xone\preceqqnu_{\og} \xtwo$ and $c(\xone)\preceqqnu_{\og} c(\xtwo)$, respectively, if and only if for \emph{every} $\vR\in\VR_\geqqV$, it holds that $\vR^\top c(\xone)\leq \vR^\top c(\xtwo)$.
	\item $\xone$ \emph{$\og$-ordinally dominates} $\xtwo$ and $c(\xone)$ \emph{$\og$-ordi\-nally dominates} $c(\xtwo)$, denoted by $\xone\preceqnu_{\og} \xtwo$ and $c(\xone)\preceqnu_{\og} c(\xtwo)$, respectively, if and only if  $\xone$ weakly $\og$-ordinally dominates $\xtwo$ and there exists $\vR^*\in\VR_\geqqV$ such that $\nu^{*\top}_{\og} c(\xone)< \nu^{*\top}_{\og} c(\xtwo)$.
	\item $x^*\in\zulM$ is called \emph{$\og$-ordinally efficient}, if there does not exist an $x\in\zulM$ such that $x\preceqnu_{\og} x^*$.
	\item $c(x^*)$ is called \emph{$\og$-ordinally non-dominated outcome vector} of problem~\eqref{eq:WOOP}, if $x^*$ is $\og$-ordinally efficient.
\end{enumerate}
\end{defi}

In the following we introduce ordering cones to reformulate the $\og$-dominance relation by a weighted ordering cone. 
Therefore, we use a characterization of multi-objective efficiency by means of weighted-sum scalarizations.

\section{Ordering Cones and Efficiency}\label{sec:cones}
Many order relations in multi-objective optimization can be represented by dominance cones, which is also the case for the (weighted) ordinal dominance relation. 
This conic representation provides deeper insights of weighted ordinal dominance and its interrelation with standard dominance concepts like, e.\,g., Pareto-dominance, see Section~\ref{subsec:dominance}.

In this section we summarize some basic results on binary relations and cones, recall dominance based on cones and analyze the ordering cone for the $\og$-dominance relation.
An in-depth introduction to polyhedral and convex cones is given, e.\,g., in \cite{rockafellar1970a} and \cite{yu74cone}.

\subsection{Binary relations and Cones}
We recall the interrelation between binary relations and cones as well as some results on (polyhedral) cones.
For a self-contained introduction to ordering cones, binary relations and their use in multi-objective optimization, see, for example, \cite{Ehrg05,engau2007domination,jahn11vector,ziegler95}.
A binary ordering relation $\mathcal{R}\subseteq\R^K\times\R^K$ is a set of pairs of vectors in $\R^K$. 
Note that we denote $(v^1,v^2)\in\mathcal{R}$ also as $v^1\mathcal{R} v^2$ for $v^1,v^2\in\R^K$ in the following.
Binary relations are closely related to cones. A \emph{cone} $C \subseteq\R^K$ is a nonempty subset of $\R^K$ such that $\lambda \, u  \in C$ for all $u\in C$ and for all $\lambda\in\R$, $\lambda>0$.

Every cone $C\in\R^K$ induces a binary (ordering) relation $\mathcal{R}_C\subseteq\R^K\times\R^K$ by $(u,v)\in\mathcal{R}_C$ if and only if $(v-u)\in C$.  The induced binary relation is compatible with scalar multiplication, i.\,e., $(u,v)\in\mathcal{R}_C$ implies $(\lambda\,  u,\lambda\, v)\in\mathcal{R}_C$ for all $u,v\in\R^K$ and $\lambda>0$, and addition, i.\,e., $(u+z,v+z)\in\mathcal{R}_C$ for all $(u,v)\in\mathcal{R}_C$ and for all $z\in \R^K$.
Conversely, every binary relation $\mathcal{R}\subseteq\R^K\times\R^K$ on $\R^K$ which is compatible with scalar multiplication and addition, i.\,e., $(u+z,v+z)\in\mathcal{R}$ for all $z\in\R^K$ and all $(u,v)\in\mathcal{R}$, induces the cone $C_\mathcal{R}\coloneqq\{(v-u)\in\R^K:(u,v)\in\mathcal{R}\}$.

Due to this interrelation, there exists a close connection between properties of binary relations and properties of cones, see, e.\,g., \cite{Ehrg05}. Let the binary relation $\mathcal{R}\subseteq\R^K\times\R^K$ be compatible with scalar multiplication and addition, than it holds that
\begin{itemize}
    \item $\mathcal{R}$ is \emph{reflexive}, i.\,e., $v\mathcal{R}v$ for all $v\in\R^K$, if and only if $0\in C_\mathcal{R}$,
    \item $\mathcal{R}$ is \emph{transitive}, i.\,e., $v^1\mathcal{R}v^2$ and $v^2\mathcal{R}v^3$ implies $v^1\mathcal{R}v^3$ for all $v^1,v^2,v^3\in\R^K$, if and only if $C_\mathcal{R}$ is \emph{convex}, i.\,e., $\lambda u^1+(1-\lambda) u^2\in C$ for all $u^1,u^2\in C$ and for all $\lambda\in[0,1]$,
    \item $\mathcal{R}$ is \emph{antisymmetric}, i.\,e., $v^1\mathcal{R}v^2$ and $v^2\mathcal{R}v^1$ implies $v^1=v^2$ for all $v^1,v^2\in\R^K$, if and only if $C_\mathcal{R}$ is \emph{pointed}, i.\,e., for $u\in C\setminus\{0\}$ it holds that $-u\notin C$.
\end{itemize}

We recall the definition of polyhedral cones as well as the characterization of their dual cones. The \emph{dual cone} of a given cone $C\subset \R^K$ is defined as $C^*\coloneqq\{d\in\R^K:d^\top c\geq0\text{ for all }c\in C\}$.

\begin{theorem}[c.f.~\citealt{rockafellar1970a,boyd2004convex}]\label{thm:coneProp}
    Let $C$ be a closed convex cone. Then it holds:
    \begin{itemize}
        \item The dual cone of the dual cone of a closed convex cone is the cone itself, i.\,e., $(C^*)^*=C$.
        \item If $\clos(C)$ is pointed, then $\interior(C^*)\neq\emptyset$.
\end{itemize}
\end{theorem}

A specific type of cones are the \emph{polyhedral cones}, which can be described as the intersection of $p$ closed halfspaces. If the normal vectors of the  corresponding defining hyperplanes are taken as rows of a matrix $A\in\R^{p\times K}$, then a polyhedral cone is defined as $\hcone(A)\coloneqq\{y\in\R^K:A\,y\geqq0\}$. Similarly, polyhedral cones can be described by their $m$ extreme rays, which form the columns of a matrix $B\in\R^{K\times m}$. Then the polyhedral cone is defined as 
$\vcone(B)\coloneqq\bigl\{B\lambda:\lambda\in\R^m,\lambda\geqq 0\bigr\}$. The Weyl-Minkowski-Theorem, cf.\ \citealt{ziegler95}, states that a cone $C\subseteq \R^K$ is finitely generated by $m$ vectors in $\R^K$ if and only if it is a finite intersection of $p$ halfspaces in $\R^{K}$. The following theorem provides an explicit description of the dual cone of a polyhedral cone and it states that the extreme rays of a cone correspond to the normal vectors of its dual cone.

\begin{theorem}[\citealt{Fukuda2016Lecture}]\label{thm:dualPolyhCone}
    Let $B\in\R^{K\times m}$ be a matrix. Then it holds $$(\vcone(B))^*=\hcone(B^\top).$$
\end{theorem}

\subsection{Dominance Based on Ordering Cones}\label{subsec:dominance}
For reflexive, transitive and antisymmetric ordering relations, which are compatible with scalar multiplication and addition, dominance can be defined by their induced cone in the following way.
For a comprehensive introduction to vector optimization and in particular to general ordering cones see, e.\,g., \cite{tammer03theory,jahn11vector}. 
\begin{defi}[c.f.~\citealt{engau2007domination,Klamroth2023Ordinal}] \label{def:coneOpt}
	Let $Y\subset\R^K$ be a nonempty set and let $C_\mathcal{R}\subset\R^K$ be a cone induced by a  partial order 
	$\mathcal{R}\subset\R^K\times \R^K$ (i.\,e., $\mathcal{R}$ is reflexive, transitive and antisymmetric) which is compatible with scalar multiplication and addition. Then the sets
	\begin{align*}
	\Nset(Y,C_\mathcal{R})&\coloneqq\{y\in Y: (y-C_\mathcal{R})\cap Y=\{y\}\} \\
	\Nset_w(Y,C_\mathcal{R})&\coloneqq\{y\in Y: (y-\interior(C_\mathcal{R}))\cap Y=\varnothing\}
	\end{align*}
	are denoted as \emph{$C_\mathcal{R}$-non-dominated} and \emph{weakly $C_\mathcal{R}$-non-dominated} set of $Y$, respectively. The corresponding pre-images $x\in\zulM$ are called \emph{$C_\mathcal{R}$-efficient} and \emph{weakly $C_\mathcal{R}$-efficient}, respectively.
Furthermore, we say that $u$ \emph{$C_\mathcal{R}$-dominates} $v$, denoted by $u\leqslant_{C_\mathcal{R}}v$, if 
	$u\in v- (C_\mathcal{R}\setminus\{0\})$, and that $u$ \emph{weakly $C_\mathcal{R}$}-dominates $v$, denoted by $u\leqq_{C_\mathcal{R}}v$, if $u\in v-C_\mathcal{R}$.
\end{defi}
 Note that a cone $C_\mathcal{R}\subset\R^K$ induced by a  partial order $\mathcal{R}\subset\R^K\times \R^K$ is always a closed and convex cone.

An important special case is Pareto-dominance, which is often used in multi-objective optimization. The Pareto cone $P$ is a polyhedral cone that is induced by the identity matrix, i.\,e., $P=\hcone(I)=\vcone(I)=\{x\in\R^K\colon x_i\geq 0 \;\forall i\in\{1,\ldots, K\}\}$. We say an outcome vector $\yone$ of the outcome set $Y\subseteq\R^K$ Pareto dominates another outcome vector $\ytwo$ ($\yone\leq_P\ytwo)$ if and only if $\yone\in \ytwo- (P\setminus\{0\})$, i.\,e., $\yone_i\leq\ytwo_i$ for all $i=1,\dots, K$ and $\yone\neq \ytwo$. Similarly, weak Pareto dominance is denoted by $\yone\leqq_P\ytwo$.

If we consider a multi-objective optimization problem with a pointed polyhedral ordering cone, then we can transform the problem into a (possibly higher dimensional) multi-objective optimization problem with the Pareto cone.
\begin{theorem}[see, e.\,g., \citealp{sawaragi1985theory},\citealp{engau2007domination}]\label{thm:linearTransformation}
	Let $Y\subset\R^K$ be a nonempty set and let $\hcone(A)$ be a pointed cone induced by a matrix $A\in\R^{p\times K}$, i.\,e., $\rank(A)=K$. Then it holds that
	\[
			A\cdot \Nset(Y,\hcone(A))= \Nset(A\cdot Y,\pareto).
	\]\end{theorem}
Note that if the matrix $A$ is invertible, then the non-dominated set of the original problem can be computed as $\Nset(Y,\hcone(A))= A^{-1}\cdot \Nset(A\cdot Y,\pareto)$.

The following results characterize cone dominance through weighting vectors in the associated dual cone. We first consider the case $\leqq_C$ (where both vectors may be equal) and then transfer the result to $C$-dominance, i.\,e., to the case $\leqslant_C$. These results are closely related to the well-known ``exactness'' and ``completeness''-properties of weighted sum scalarizations that can be found in several sources, see, e.\,g., \cite{Dempe2015effects,Ehrg05,jahn11vector,miettinen1998}
and the references therein.
\begin{theorem}
\label{thm:weighted}
    Let $\yone,\ytwo\in\R^K$ and let \(C\subset\R^K\) be a closed, convex cone. Then $\yone$ weakly \(C\)-dominates \(\ytwo\), $\yone\leqq_C \ytwo$ (i.\,e., \(\yone\in \ytwo - C\)) if and only if $\lambda^\top\yone\leq\lambda^\top\ytwo$ for all weighting vectors $\lambda\in C^*$.\end{theorem}

\begin{proof}
	First we show that $\yone\leqq_C \ytwo$ implies that $\lambda^\top \yone\leq\lambda^\top \ytwo$ for all $\lambda\in C^*=\{\lambda\in\R^K:\lambda^\top y\geq 0\text{ for all }y\in C\setminus\{0\}\}$. Indeed, if $\yone\leqq_C \ytwo$, then by definition there exists \(c\in C\) such that \(\yone=\ytwo -c \). If $c=0$ then it holds that $\lambda^\top \yone = \lambda^\top \ytwo$ for all $ \lambda\in C^*$. If $c\neq 0$ then we have
    $\lambda^\top \yone = \lambda^\top \ytwo - \lambda^\top c$ and, since $\lambda^\top c\geq 0$ for all $\lambda\in C^*$, it follows that $\lambda^\top \yone \leq \lambda^\top \ytwo$. 

	To show the other direction we assume, to the contrary, that there exist $\yone,\ytwo\in \R^K$ such that $\lambda^\top \yone \leq \lambda^\top \ytwo$ for all $\lambda\in C^*$ but $\yone \nleqq_C  \ytwo$. From $\yone \nleqq_C  \ytwo$ it follows that
    \( c\coloneqq \ytwo - \yone \notin C\). 
    Then $\lambda^\top \yone \leq \lambda^\top \ytwo$ for all $\lambda\in C^*$ implies that 
    \[
      \lambda^\top c = \lambda^\top \ytwo - \lambda^\top \yone \geq 0 \quad \forall \lambda\in C^*. 
    \]
    Consequently, \(c\) is an element of \((C^*)^*\), the dual cone of \(C^*\). Since \(C\) is closed and convex, Theorem~\ref{thm:coneProp} implies that $(C^*)^*=C$. Hence, \(c\in C\) which is a contradiction.
\end{proof}

\begin{lem}\label{lem:weighted2}
    For a closed, convex and pointed cone $C$,
    a point $\yone$ $C$-dominates $\ytwo$, i.\,e., $\yone\leqslant_C\ytwo$ if and only if $\lambda^\top\yone\leq\lambda^\top\ytwo$ holds for all $\lambda\in C^*$ and there exists a weighting vector $\bar{\lambda}\in C^*$ such that $\bar{\lambda}^\top\yone<\bar{\lambda}^\top\ytwo$.
\end{lem}

\begin{proof}
    First note that $\interior(C^*)\neq\emptyset$ since $C$ is pointed, cf.\ Theorem~\ref{thm:coneProp}. Then
    the result follows from Theorem~\ref{thm:weighted} and the fact that $\yone\leqslant_C \ytwo$ if and only if $\yone\leqq_C \ytwo$ and $\yone\neq \ytwo$. Indeed, $\yone\leqslant_C \ytwo$ implies that $\yone=\ytwo-c$ with $c\neq 0$ and hence $\lambda^\top c\geq 0$ for all $\lambda\in C^*$ and $\lambda^\top c > 0$ for all $\lambda\in \interior(C^*)$. Hence, 
    $\bar{\lambda}^\top \yone<\bar{\lambda}^\top \ytwo$ for all $\bar{\lambda}\in\interior(C^*)\neq\emptyset$ in this case. 

    Conversely, if $\lambda^\top \yone\leq\lambda^\top \ytwo$ for all $\lambda\in C^*$ and there exists  $\bar{\lambda}\in C^*$ such that $\bar{\lambda}^\top \yone<\bar{\lambda}^\top \ytwo$, then $\yone\neq \ytwo$ and $\yone\leqq_C \ytwo$ by Theorem~\ref{thm:weighted},
which concludes the proof.
\end{proof}

\section{The Weighted Ordinal Ordering Cone}\label{sec:weightedCone}
Starting with the assumption $\og\in\bar{\Omega}$, we later restrict the choice of weights to $\omega_i\cdot \gamma_i <1$ for all $i=1,\ldots,K-1$ in order to obtain a pointed ordering cone.
We investigate general properties of the relation $\preceqqnu_{\og}$ and the corresponding weighted ordinal ordering cone $\tailC$. In particular, we derive the description of this cone by its facets. Moreover, we investigate  specific choices of the weights and the interrelation of the ordinal ordering relation with other well-known ordering relations.

\subsection{Properties} 
According to Theorem~\ref{thm:weighted}, the weighted ordinal dominance relations (Definition~\ref{def:ordOpt-wg}) are defined by modifications of their dual cone, as compared to the standard ordinal dominance (Definition~\ref{def:ordOpt}). Note that the concept of incorporating preference information in multi-objective programming problems by modifying the dual of the ordering cone is investigated in detail in  \citet{Kaddani2017Weighted}.

The additional constraints (\ref{eq:3}--\ref{eq:4}) for the set of weighted numerical representations in Definition~\ref{def:generalNumRep}  correspond to linear constraints of the dual of the weighted ordinal dominance cone, i.\,e.,  $\VR_\geqqV=\tailC^*$. 
Moreover, one can show using Theorem~\ref{thm:dualPolyhCone}, that the cone $\tailC$ is spanned by the vectors $u^i\in\R^K$ and $g^i\in\R^K$ for all $i=1,\dots,K-1$, as defined in equations \eqref{eq:u} and \eqref{eq:g}. Thus, it holds that:
\[\tailC=\vcone(B)  \quad\text{with}\; B=(u^1,\ldots,u^{K-1},g^1,\ldots,g^{K-1})\in\smash{\R^{K\times 2\,(K-1)}}.\]

Note that some of the spanning vectors in the matrix $B$ may not be extreme rays. Thus, some of the columns of $B$ may be redundant and can be removed.
In particular, if $\g_i=0$ for some $i\in\{2,...,K-1\}$, then the vectors $u^{i-1},g^{i-1},g^i$ are linearly dependent. Analogously, one can show that whenever a weight $\w_i=0$ for some $i\in\{1,...,K-2\}$, then the vector $u^i$ lies in the hyperplane which contains $u^{i+1}$ and $g^{i+1}$.

\begin{lem}\label{lem:pointedCone}
    Let $\og\in\bar{\Omega}$.
    The ordering cone $\tailC$ is pointed (that is, the relation $\preceqqnu_{\og}$ is antisymmetric) if and only if $\w_i\,\g_i<1$ for all $i=1,\dots,K-1$.
\end{lem}
\begin{proof}
  We first show that $\w_i\,\g_i<1$ for all $i=1,\dots,K-1$ implies that the ordering cone $\tailC$ is pointed. 
  We consider the following matrix 
  \begin{equation*}
		M=\begin{pmatrix}
		    m_1^\top\\ \vdots\\m_K^\top
		\end{pmatrix} \coloneqq
		\begin{tikzpicture}[baseline=(current bounding box.center)]
		\matrix (m) [matrix of math nodes,nodes in empty cells,right delimiter={)},left delimiter={(},inner sep=-2.5pt,nodes={inner sep=1.25ex}]{
			1 & \w_1 & \w_1\,\w_2&  & \prod_{\ell=1}^{K-1}\w_\ell\\
			\g_1& 1 &\w_2& & \prod_{\ell=2}^{K-1}\w_\ell \\
			&  & & &\\
			& & & 1& \w_{K-1}\\
			\prod_{\ell=1}^{K-1}\g_\ell& & & \g_{K-1}&1 \\
		} ;
		\draw[loosely dotted,thick] (m-2-2)-- (m-4-4);
		\draw[loosely dotted,thick] (m-2-3)-- (m-4-5);
		\draw[loosely dotted,thick] (m-1-3)-- (m-1-5);
		\draw[loosely dotted,thick] (m-2-3)-- (m-2-5);
		\draw[loosely dotted,thick] (m-2-5)-- (m-4-5);
		\draw[loosely dotted,thick] (m-2-1)-- (m-5-4);
		\draw[loosely dotted,thick] (m-2-1)-- (m-5-1);
		\draw[loosely dotted,thick] (m-5-1)-- (m-5-4);
		\end{tikzpicture},
		\end{equation*}
    According to Definition~\ref{def:generalNumRep}, all row vectors \(m_i^\top\) of $M$ are contained in $\VR_\geqqV$. Hence, it holds that $\vcone(M)\subseteq \VR_\geqqV$. With Theorem~\ref{thm:coneProp} we obtain $\bigl(\VR_\geqqV\bigr)^*\subseteq \bigl(\vcone(M)\bigr)^*$, which is equivalent to $\tailC\subseteq \hcone(M^\top)$.  

    Next, we show that $\hcone(M^\top)$ is pointed, which implies that $\tailC$ is pointed, as well. Let $v^1,v^2\in\R^K$ and assume that $v^1\preceqqnu_{\og} v^2$ and $v^2 \preceqqnu_{\og} v^1$. 
    Since $\tailC \subseteq \hcone(M^\top)$ it holds that $v^2 - v^1\in \hcone(M^\top)$ and $v^1 - v^2\in \hcone(M^\top)$ and consequently $M^\top(v^2-v^1)=0$. 
    The matrix $M$ is quadratic and has full rank since $\w_i<\frac{1}{\g_i}$ for all $i=1,\dots,K$.  This implies $v^1=v^2$ and thus the relation induced by the cone $\hcone(M^\top)$ is antisymmetric, i.\,e., $\hcone(M^\top)$ is pointed.

    The other direction of the statement is proven by contradiction. We assume that there exists an index $j\in\{1,\dots,K-1\}$ such that $\w_j\,\g_j=1$, i.\,e., $\w_j\, \nu_{j}=\nu_{j+1}$. Then, it holds that $u^j=-\w_j\,g^j$ with $u^j,\,g^j$ as defined in \eqref{eq:u} and \eqref{eq:g}. Thus, we can deduce that $u^j,-u^j\in\tailC$ with $u^j\neq 0$. This implies that the cone $\tailC$ is not pointed, which completes the proof.
\end{proof}

\begin{rem}  
    In the case of $\w_i=1/\g_i$ for some $i\in\{1,\dots,K-1\}$ the corresponding dominance cone contains a line through the origin and is therefore not pointed. In this situation we can merge categories (as described in Remark~\ref{rem:merge}) and obtain a lower-dimensional problem with a pointed dominance cone.
\end{rem} 

Thus, we assume in the following that the cone $\tailC$ is pointed  and let $A$ be the matrix such that $\tailC=\hcone(A)$, $A\in\R^{p\times K}$, i.\,e., the matrix that describes the weighted ordinal cone by its $p$ facets.
By applying Theorem~\ref{thm:linearTransformation}, we can compute the non-dominated set of problem~\eqref{eq:WOOP} by transforming the problem into a potentially higher dimensional problem w.\,r.\,t.\ Pareto dominance, i.\,e.,
    \[
		A\cdot \Nset(Y,\hcone(A))= \Nset(A\cdot Y,\pareto),
	\]
    where \(P\) denotes the Pareto cone.

The weighted ordering cone $\tailC$ has up to $2(K-1)$ extreme rays. Thus, for $K\geq 2$ there exists in general no bijective linear transformation of \eqref{eq:WOOP} to a Pareto optimization problem as in the case of standard ordinal optimality \citep[see][]{Klamroth2023Ordinal}. The resulting problem w.\,r.\,t.\ Pareto optimality is in general higher-dimensional and depends on the number of facets of the weighted ordinal ordering cone $\tailC$.

\subsection{Description of the Weighted Ordinal Ordering Cone by its Facets}\label{subsec:descr}
In the following, we derive an alternative description of the weighted ordinal ordering cone $\tailC=\{(v-u)\in\R^K:u\preceqqnu_{\og}v\}$ based on its facets, i.\,e., we determine a matrix \(A\) such that \(\tailC=\hcone(A)\).
In general, this requires the double description method; see, e.\,g., \cite{DoubleDescriptionRevisited}.

Recall that the dual of the weighted ordinal ordering cone $\tailC$ represents the set of numerical representations. We can interpret the weights $\omega,\ \gamma$ as so-called ratio bounds between the weights of consecutive categories, i.\,e., 
\[
  \g_i\leq\frac{\nu_i}{\nu_{i+1}}\leq \frac{1}{\w_i}\ \quad \forall i=1,\dots,n-1.
\]
Ratio bounds are investigated in detail in \cite{Ahn2014Establishing}, where the following result is proven:
\begin{lem}[\cite{Ahn2014Establishing}] \label{lem:ahn}
    Let the weights $\omega,\ \gamma$ be strictly positive and $\w_i\,\g_i<1$ for $i=1,\dots,K-1$, i.\,e., $\og\in\{\og\in\Omega:0<\w_i,\ 0<\g_i\text{ for all }i=1,\dots,K-1\}$. Then the number of extreme points of the polytope $\mathcal{P}=\{\nu\in\R^K:\g_i\leq\frac{\nu_i}{\nu_{i+1}}\leq \frac{1}{\w_i}\ \forall i=1,\dots,n-1,\ ||\nu||_1=1\}$ is equal to $2^{K-1}$, and each extreme point can be determined by solving the following system of linear equations:
    $$\frac{\nu_i}{\nu_{i+1}} = \frac{1}{\w_i} \;\vee\; \frac{\nu_i}{\nu_{i+1}} = \g_i \quad 
     ,\ 1\leq i \leq K-1,\ \sum_{i=1}^K \nu_i=1.$$
\end{lem}
Obviously, the polytope $\mathcal{P}$ from Lemma~\ref{lem:ahn} has at most $2^{K-1}$ extreme points, as either $\frac{\nu_i}{\nu_{i+1}}=\g_i$ or (exclusively) $\frac{\nu_i}{\nu_{i+1}}=\frac{1}{\w_i}$ for $i=1,\dots,K-1$, due to the constraint $\w_i\,\g_i<1$. In the following, we describe the $2^{K-1}$ extreme points explicitly. Note that the extreme points of the polytope $\mathcal{P}$ correspond to the extreme rays of $\VR_\geqq$ and the normal vectors of the facets of $\tailC$. As we are interested in the normal directions of the facets of $\tailC$, we can ignore the lengths of the vectors and omit the equation $\sum_{i=1}^K \nu_i=1$. Thus, we obtain the following description of the normal vectors.

\begin{lem}\label{thm:normalvectors}
    Let a weighted ordinal ordering cone $\tailC$ be given for $\og\in\Omega$. Let $r^i\in\{u^i,g^i\}$, $i=1,\dots,K-1$, be a subset of \(K-1\) extreme rays defining a facet \(\mathcal{F}\), see equations~\eqref{eq:u} and \eqref{eq:g}. 
    Then the normal vector $n=(n_1,\ldots,n_K)^\top$ of \(\mathcal{F}\) is given by $n_k=\prod_{i=1}^{K-1} d^i_k$ with vectors $d^i\in\R^K$ that are defined depending on the corresponding extreme rays $r^i$ 
    for $i=1,\dots,K-1$, i.\,e., 
    \[    
        d^i_k=\begin{cases*}
            \w_i & if  $r^i=u^i$  and  $k> i$\\
            1 &  if   $r^i=u^i$  and $k \leq i$\\
            \g_i & if  $r^i=g^i$  and  $k\leq i$\\
            1 & if  $r^i=g^i$  and $k > i$ .
        \end{cases*}
    \]
    
\end{lem}
\begin{proof}
    We will show that \(n\) is orthogonal to all extreme rays incident to \(\mathcal{H}\), i.\,e., $n^\top r^i=0$ holds for all $i=1,\dots,K-1$. This can be verified by a straightforward calculation.
    
    Every extreme ray \(r^i\) has at least one and at most two non-zero entries, namely in the \(i\)-th and/or \((i+1)\)-st component. The corresponding components of the normal vector \(n\), differ in exactly one factor of the product.
For $r^i=u^i$ it holds that $n_{i+1}=\w_i\,n_i$ and thus $n^\top u^i = -\w_i\, n_i + \,n_{i+1}= 0$. For $r^i=g^i$ it holds that $n_{i}= \g_i\,n_{i+1}$ and thus $n^\top g^i= \,n_i-\g_i\, n_{i+1}=0$.
\end{proof}

\begin{example}
Let a weighted ordinal ordering cone $\tailC$ for $K=4$ and $\og\in\Omega$. Then, there is one facet described by the extreme rays $g^1,u^2,g^3$. The corresponding normal vector of the facet is defined by $$d^1=\begin{pmatrix}
    \g_1\\ 1 \\ 1 \\1 
\end{pmatrix},\ d^2=\begin{pmatrix}
    1\\ 1\\ \w_2\\ \w_2
\end{pmatrix},\ d^3=\begin{pmatrix}
    \g_3 \\ \g_3 \\ \g_3\\ 1 
\end{pmatrix},\text{ and thus }n=\begin{pmatrix}
    \g_1\,\g_3\\ \g_3 \\ \w_2\,\g_3 \\ \w_2 
\end{pmatrix},$$ as it holds that 
    \begin{equation*}
        \begin{pmatrix}
            1 & -\g_1 & 0 & 0 \\
            0 & -\w_2 & 1 & 0 \\
            0 & 0 & 1 & -\g_3 
        \end{pmatrix}\cdot
        \begin{pmatrix}
            \g_1\,\g_3\\
            \g_3\\
            \w_2\,\g_3\\
            \w_2
        \end{pmatrix}=
        \begin{pmatrix}
            0\\0\\0\\0
        \end{pmatrix}.
    \end{equation*}
    Similarly, we get 
    the weighted ordinal ordering cone $\tailC=\hcone(A)$ for $K=4$, $\og\in\Omega$ with strictly positive weights, i.\,e., $\w_i>0,\ \g_i>0$ for all $i=1,\dots,K-1$, defined by
    $$A\coloneqq\begin{pmatrix}
        1 & \w_1 & \w_1\,\w_2 & \w_1\,\w_2\,\w_3\\
        \g_3 & \g_3\,\w_1 & \g_3\,\w_1\,\w_2 & \w_1\,\w_2\\
        \g_2 & \g_2\,\w_1 & \w_1 & \w_1\,\w_3\\
        \g_2\,\g_3 & \g_2\,\g_3\,\w_1 & \g_3\,\w_1 & \w_1\\
        \g_1 & 1 & \w_2 & \w_2\,\w_3\\
        \g_1\,\g_3 & \g_3 & \g_3\,\w_2 & \w_2\\
        \g_1\,\g_2 & \g_2 & 1 & \w_3\\
        \g_1\,\g_2\,\g_3 & \g_2\,\g_3 & \g_3 & 1
    \end{pmatrix}.$$
\end{example}

Note that, if some of the weights \(\g_i,\w_i\) \(i\in\{1,\dots,K-1\}\) are zero, the total number of facets decreases.  In the following theorem the number of facets is determined for the case that $\w_i>0$ for all $i=1,\dots,K-1$, which holds, for example, for strictly ordered categories (i.\,e., $\w_i>1$ for all $i=1,\dots,K-1$).

\begin{theorem}\label{thm:numberfacets0} 
    Let $\og\in\Omega$ and let
    $J=\{j_1,\dots,j_\ell\}\subseteq \{1,\dots,K-1\}$ be the set of indices with $\g_j=0$ for all $j\in J$ and $j_k<j_{k+1}$ for all $k=1,\dots,\ell-1$. Moreover, let $\g_i>0$ for all $i\in\{1,\dots,K-1\}\setminus J$ as well as $\w_i>0$ for all $i=1,\dots,K-1$. Then the weighted ordinal ordering cone $\tailC$ has 
    \begin{equation*}
        2^{K-1-\ell} + \sum_{k=1}^\ell 2^{K-1-j_k-(\ell-k)}
    \end{equation*}
    facets.
\end{theorem}

\begin{proof}
    We consider all possible combinations of (not necessarily extreme) rays \(g^i, u^i\), \(i\in\{1,\ldots,K-1\}\), which may describe a facet.
    First we choose $r^j=u^j$ for all $j\in J$. Thus, there are $K-1-\ell$ extreme rays $r^i$ left, which can be either $u^i$ or $g^i$, which leads to $2^{K-1-\ell}$ possible combinations of extreme rays. 
    
    For a fixed index $k\in\{1,\dots,\ell\}$ we consider the facets with spanning ray $g^{j_k}$. The construction of the normal vectors in Lemma~\ref{thm:normalvectors} implies that the components \(n_i=0\) for all \(i=1,\ldots, j_k\). For all indices \(i<j_k\) the selection of the spanning ray \(u^i\) instead of \(g^i\) yields an stretched normal vector (by the factor of \(\w_{i}\)), representing the same facet. Consequently, we can fix $r^i=g^i$ for $i=1,\dots,j_k$. 
    With the same argument, we can fix all \(r^j=u^j\) for \(j\in J\) and \(j>j_k\).
    Thus, only the selection of $r^i\in\{u^i,g^i\}$ if $i>j_k$ and $i\notin J$ yields different facets, which are in total $K-1-j_k-(\ell-k)$ rays. Hence, we get $2^{K-1-j_k-(\ell-k)}$ additional facets. We get the final formula by summing up over all possible choices of $k$ and adding the number of combinations with $r^j=u^j$ for $j\in J$.
\end{proof}

Note that if $\g_i=0$ for all $i=1,\dots,K-1$ it holds that  $\ell=K-1$ and $j_k=k$ for $k=1,\dots,\ell$ and thus the cone has 
\begin{equation*}
        2^{K-1-(K-1)} + \sum_{k=1}^{K-1} 2^{K-1-k-(K-1-k)}=2^0+ \sum_{k=1}^{K-1} 2^0=K
\end{equation*}
facets.
\begin{example}
    We consider a weighted ordinal ordering cone $\tailC$ for $K=5$ with $\g_1=0$, $\g_3=0$ and all other weights are strictly positive. Then the corresponding matrix $A$ is given by 
    \begin{equation*}
        A=\begin{pmatrix}
            1 & \w_1 & \w_1\,\w_2 & \w_1\,\w_2\,\w_3 & \w_1\,\w_2\,\w_3\,\w_4\\
            \g_2 & \g_2\,\w_1 & \w_1 & \w_1\,\w_3 & \w_1\,\w_3\,\w_4\\
            \g_4 & \g_4\,\w_1 & \g_4\,\w_1\,\w_2 & \g_4\,\w_1\,\w_2\,\w_3 & \w_1\,\w_2\,\w_3 \\
            \g_2\,\g_4 & \g_2\,\g_4\,\w_1 & \g_4\,\w_1 & \g_4\,\w_1\,\w_3 & \w_1\,\w_3\\
            0 & 1 & \w_2 & \w_2\,\w_3 & \w_2\,\w_3\,\w_4 \\
            0 & \g_2 & 1 & \w_3 & \w_3\,\w_4 \\
            0 & \g_4 & \g_4\,\w_2 & \g_4\,\w_2\,\w_3 & \w_2\,\w_3 \\
            0 & \g_2\,\g_4 & \g_4 & \g_4\,\w_3 & \w_3 \\
            0 & 0 & 0 & 1 & \w_4\\
            0 & 0 & 0 & \g_4 & \w_1\\
        \end{pmatrix}
    \end{equation*}
    The corresponding ordinal cone has $$2^{4-2}+2^{4-1-(2-1)}+2^{4-3-(2-2)}=2^2+2^2+2^1=10$$
    facets.
\end{example}

\subsection{Special Cases}\label{subsec:specialCases}
In general, the weighted ordinal ordering cone $\tailC=\vcone(M)$ has up to $2^{K-1}$ facets and hence the associated multi-objective optimization problem has as many objective functions. Nonetheless, there exist some specific cases for which the number of facets equals $K$, i.\,e., the dimension of the transformed multi-objective problem w.\,r.\,t.\ the Pareto cone equals the number of categories. We investigate the following selected special cases:
\begin{itemize}
    \item If $K=2$ the transformation matrix is $$ A=\begin{pmatrix}
        1 & \w_1\\ \g_1 & 1
    \end{pmatrix}.$$ Note that $p=2$ is a special case as all polyhedral cones, including the Pareto cone, can be modeled by appropriately weighted ordinal cones. This does not hold for higher dimensions.
    The result of the special case is aligned with the two-dimensional special case in \cite{Hunt2010Relative}, where allowable trade-offs between objectives are considered, which (in contrast to our approach) measure the deterioration of one objective (or a group of objective) in relation to the improvement of all other objectives. However, in the case of only two objectives this coincides with our approach.
    \item If $\w=0$ the transformation matrix is  $A^{(\zeros,\g)}=(a_{ij})_{i,j=1,\dots,K}\in\R^{K\times K}$ with \begin{equation*}a_{ij}=\begin{cases*}
		0,& if $i< j$\\
		1,& if $i=j$\\
		\prod_{\ell=j}^{i-1}\g_\ell,& if $i> j$
		\end{cases*},
		\end{equation*}
		\begin{equation*}\text{i.\,e., }
		A^{(\zeros,\g)}=
		\begin{tikzpicture}[baseline=(current bounding box.center)]
		\matrix (m) [matrix of math nodes,nodes in empty cells,right delimiter={)},left delimiter={(},inner sep=-2pt, column sep=3pt, nodes={inner sep=1.5ex}]{
			1 & 0 & &  & 0\\
			\g_1&  && &  \\
			&  & & &\\
			& & & & 0\\
			\prod_{\ell=1}^{K-1}\g_\ell& & & \g_{K-1}&1 \\
		} ;
		\draw[loosely dotted,thick] (m-1-1)-- (m-5-5);
		\draw[loosely dotted,thick] (m-1-2)-- (m-4-5);
		\draw[loosely dotted,thick] (m-1-2)-- (m-1-5);
\draw[loosely dotted,thick] (m-1-5)-- (m-4-5);
		\draw[loosely dotted,thick] (m-2-1)-- (m-5-4);
		\draw[loosely dotted,thick] (m-2-1)-- (m-5-1);
		\draw[loosely dotted,thick] (m-5-1)-- (m-5-4);
		\end{tikzpicture}.
		\end{equation*}
        Note that by considering only the weights $\g$ we loose the order of the categories. Even more so, if $\g\geq1$ the order of the categories is reversed, i.\,e., the last category gets the best one.
        
    \item If $\g=0$ the transformation matrix is $A^{(\w,\zeros)}=(a_{ij})_{i,j=1,\dots,K}\in\R^{K\times K}$ with \begin{equation*}a_{ij}=\begin{cases*}
		\prod_{\ell=i}^{j-1}\w_\ell,& if $i< j$\\
		1,& if $i=j$\\
		0,& if $i> j$
		\end{cases*},
		\end{equation*}
		\begin{equation*}\text{i.\,e., }
		A^{(\w,\zeros)}=
		\begin{tikzpicture}[baseline=(current bounding box.center)]
		\matrix (m) [matrix of math nodes,nodes in empty cells,right delimiter={)},left delimiter={(},inner sep=-2pt, column sep=3pt, nodes={inner sep=1.5ex}]{
			\;1\; & \w_1 & \w_1\,\w_2&  & \prod_{\ell=1}^{K-1}\w_\ell\\
			0&  & \w_2& & \prod_{\ell=2}^{K-1}\w_\ell \\
			&  & & &\\
			& & & & \w_{K-1}\\
			0& & & 0&\;1\; \\
		} ;
		\draw[loosely dotted,thick] (m-1-1)-- (m-5-5);
		\draw[loosely dotted,thick] (m-2-3)-- (m-4-5);
		\draw[loosely dotted,thick] (m-1-3)-- (m-1-5);
		\draw[loosely dotted,thick] (m-2-3)-- (m-2-5);
		\draw[loosely dotted,thick] (m-2-5)-- (m-4-5);
		\draw[loosely dotted,thick] (m-2-1)-- (m-5-4);
		\draw[loosely dotted,thick] (m-2-1)-- (m-5-1);
		\draw[loosely dotted,thick] (m-5-1)-- (m-5-4);
		\end{tikzpicture}.
\end{equation*}
        This case can be relevant when considering hazardous paths. Consider, for example, a chemical accident with dangerous smoke and we need to compute a path from one place to another. Then, edges are assigned to categories according to their distance to the area of the accident, maybe under consideration of additional information like wind direction. In this case, one would prefer to use a longer detour to avoid dangerous areas. E.\,g., we would prefer \(\w_1\) green edges over one red edge.
\end{itemize}

\begin{figure}
	\centering
    \subcaptionbox{$\w=0,\, \g=0$ (Pareto Cone) \label{fig:o0g0}}	[.45\linewidth]{\begin{tikzpicture}[scale=0.65]
\fill[fill=WtalUniGruen!20,fill opacity=0.8] (3,0) -- (3,3) -- (0,3) -- (0,0)  -- cycle;
	
\draw[->] (-3.2,0) -- (3.2,0) node[right] {};
	\draw[->] (0,-1.3) -- (0,3.2) node[above] {};
	
\fill[pattern color=red,pattern=north west lines] (0,3) -- (3,3) -- (3,0) -- (0,0)  -- cycle;
	\draw[red] (0,0) -- (0,3);
	\draw[red] (0,0) -- (3,0);
\end{tikzpicture} }
    \subcaptionbox{$\w=1,\, \g=0$ (Standard Ordinal Cone) \label{fig:o1g0}}	[.45\linewidth]{\begin{tikzpicture}[scale=0.65]
\draw[->] (-3.2,0) -- (3.2,0) node[right] {};
	\draw[->] (0,-1.3) -- (0,3.2) node[above] {};
	
\fill[fill=WtalUniGruen!20,fill opacity=0.8] (3,0) -- (3,3) -- (-3,3) -- (0,0)  -- cycle;
        \draw[WtalUniGruen] (0,0) -- (3,0);
        \draw[WtalUniGruen] (0,0) -- (-3,3);
	
\fill[pattern color=red,pattern=north west lines] (0,3) -- (3,3) -- (0,0)  -- cycle;
	\draw[red] (0,0) -- (0,3);
	\draw[red] (0,0) -- (3,3);
\end{tikzpicture} }\\
	\subcaptionbox{$\w=0,\, \g=\frac{1}{4}$ \label{fig:o0g4}}[.45\linewidth]{\begin{tikzpicture}[scale=0.65]
\draw[->] (-3.2,0) -- (3.2,0) node[right] {};
	\draw[->] (0,-1.3) -- (0,3.2) node[above] {};
 
\fill[fill=WtalUniGruen!20,fill opacity=0.8] (3,-0.75) -- (3,3) -- (0,3) -- (0,0)  -- cycle;
        \draw[WtalUniGruen] (0,0) -- (3,-0.75);
        \draw[WtalUniGruen] (0,0) -- (0,3);
	
\fill[pattern color=red,pattern=north west lines] (3,0) -- (3,3) -- (0.75,3) -- (0,0)  -- cycle;
	\draw[red] (0,0) -- (0.75,3);
	\draw[red] (0,0) -- (3,0);

\end{tikzpicture} }
	\subcaptionbox{$\w=2,\, \g=0$\label{fig:o2g0}}[.45\linewidth]{\begin{tikzpicture}[scale=0.65]
\draw[->] (-3.2,0) -- (3.2,0) node[right] {};
	\draw[->] (0,-1.3) -- (0,3.2) node[above] {};
	
\fill[fill=WtalUniGruen!20,fill opacity=0.8] (3,0) -- (3,3) -- (-3,3) -- (-3,1.5)-- (0,0)  -- cycle;
        \draw[WtalUniGruen] (0,0) -- (3,0);
        \draw[WtalUniGruen] (0,0) -- (-3,1.5);
	
\fill[pattern color=red,pattern=north west lines] (0,3) -- (3,3) -- (1.5,3) -- (0,0)  -- cycle;
	\draw[red] (0,0) -- (1.5,3);
	\draw[red] (0,0) -- (0,3);
\end{tikzpicture} }\\
	\subcaptionbox{$\w=2,\, \g=\frac{1}{4}$\label{fig:o2g4}}[.45\linewidth]{\begin{tikzpicture}[scale=0.65]
\fill[fill=WtalUniGruen!20,fill opacity=0.8] (3,-0.75) -- (3,3) -- (-3,3) -- (-3,1.5) -- (0,0)  -- cycle;
        \draw[WtalUniGruen] (0,0) -- (3,-0.75);
        \draw[WtalUniGruen] (0,0) -- (-3,1.5);
	
\fill[pattern color=red,pattern=north west lines] (0.75,3) --  (1.5,3) -- (0,0)  -- cycle;
	\draw[red] (0,0) -- (0.75,3);
	\draw[red] (0,0) -- (1.5,3);
	
\draw[->] (-3.2,0) -- (3.2,0) node[right] {};
	\draw[->] (0,-1.3) -- (0,3.2) node[above] {};
\end{tikzpicture} }
    \subcaptionbox{$\w=2,\, \g=\frac{1}{2}$ (weighted sum) \label{fig:o2g2}}[.45\linewidth]{\begin{tikzpicture}[scale=0.65]
\fill[fill=WtalUniGruen!20,fill opacity=0.8] (3,-1.5) -- (3,3) -- (-3,3) -- (-3,1.5) -- cycle;
	\draw[WtalUniGruen] (0,0) -- (3,-1.5);
        \draw[WtalUniGruen] (0,0) -- (-3,1.5);
 
\draw[red] (0,0) -- (1.5,3);
	
\draw[->] (-3.2,0) -- (3.2,0) node[right] {};
	\draw[->] (0,-1.3) -- (0,3.2) node[above] {};
\end{tikzpicture} }
	\caption{The weighted ordinal ordering cones $\protect\tailC$ (green) and their corresponding dual cones (red, stripped) for different choices of $(\w,\g)$. Note that $p=2$ is a special case since all polyhedral cones can be equivalently described as weighted ordinal cones, which does not hold for higher dimensions.
	}
	\label{fig:cones}
\end{figure}
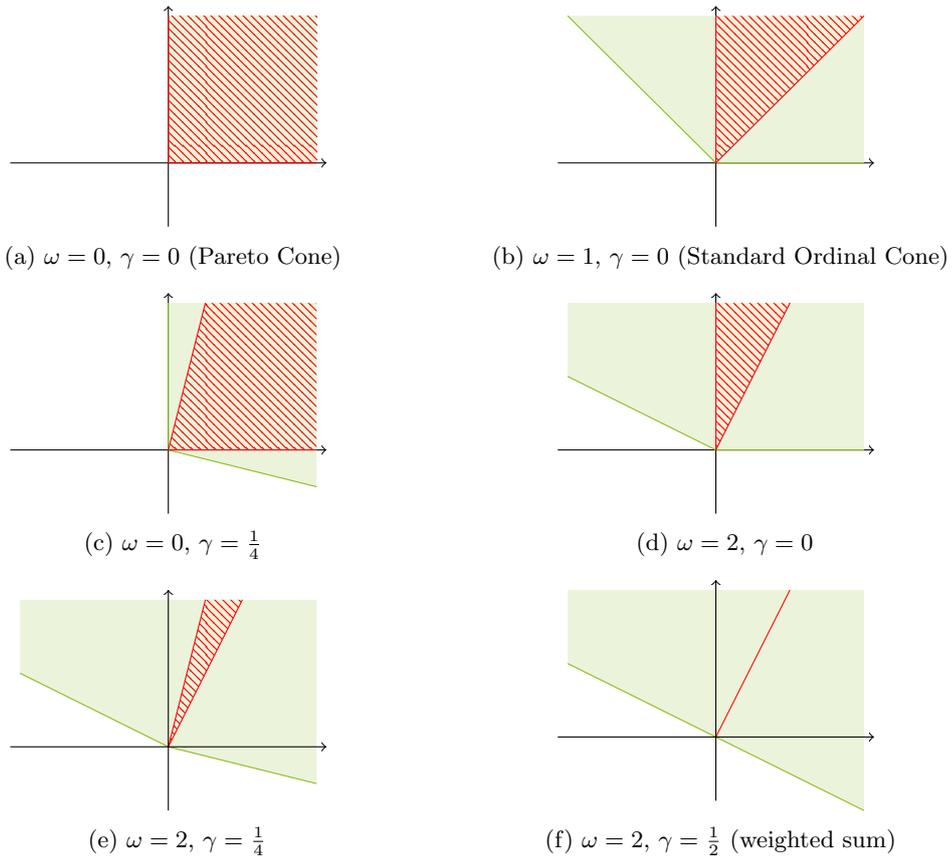

In Figure~\ref{fig:cones} the weighted ordinal ordering cone $\tailC$ and its dual cone are illustrated for $K=2$ categories and different choices of the weights $\g$ and $\w$. 
Note that for specific choices of $\w$ and $\g$ the weighted ordinal cone equals different standard cones (like the Pareto cone or the standard ordinal cone, for example). This relation holds  for an arbitrary number of categories:
    \begin{itemize}
        \item For $\w=\g=\zeros$ we get the Pareto cone $P$.
        \item For $\w=\ones$ and $\g=\zeros$ we get the standard ordinal cone.
        \item For $\w\geq \ones$ the resulting cone includes the standard ordinal cone.
        \item If  $\w_i=\frac{1}{\g_i}$ for all $i=1,\dots,K$, then the induced cone corresponds to a halfspace. This is equivalent to applying the weighted sum method with objective function $c_1(x)+\sum_{i=2}^K w_{i-1} \,c_i(x)$, $x\in \zulM$.
        \item For $\g=\zeros$ and $\w\to\infty$ the weighted ordering cone converges to the lexicographic ordering cone. The lexicographic order $\yone<_{\lex}\ytwo$ for $\yone,\ytwo\in\R^K$ holds if $\yone_j< \ytwo_j$ for some index $j\in\{1,\dots,K\}$ and $\yone_i=\ytwo_i$ for all $i=1,\dots,j-1$.
    \end{itemize}
By modifying the weights $\omega,\ \gamma$, the resulting weighted ordinal ordering cone $\tailC$ varies between the Pareto cone and the ''weighted sum cone'' (or the closure of the lexicographic ordering cone in the limit case). This interrelation of weight sets and ordering cones is investigated in a more general setting in \cite{Kaddani2017Weighted}, where the weighted sum scalarization is generalized to overcome its disadvantage that, in general, only one supported solution is computed.

\section{Numerical Results}\label{sec:numRes} 
To visualize the effect of the weights, we apply different weighted ordering cones to compute the safest paths for cyclists on two different routes. For this purpose, we use data from OpenStreetMap (OSM) \citep{OpenStreetMap} and apply an adapted version of the ordinal routing code from \cite{Santos2024Computing}. We use the python library OSMnx, see \cite{osmnx}, to generate a graph using the data from \cite{OpenStreetMap} and choose the network type ``bike'' such that paths, which cannot be used by bike, are not included in the graph. Moreover, we use the option from OSMnx to simplify the graph, see \cite{osmnx}, and we delete multiple edges between two nodes such that a simple graph is obtained. For more details on the construction of the graph, see \cite{Santos2024Computing}. 

We assign four categories to the streets, as it is done in \cite{Santos2024Computing}.
The best category contains streets with a separate bike lane, the second category contains streets with a bike lane on the street, in the third category are streets with little traffic, like e.\,g., residential streets, and in the worst category are all other streets. For the computation of efficient shortest paths, we use the multi-objective Dijkstra algorithm presented in \cite{Casas2021Improved}. 

\begin{figure}
	\centering
    \subcaptionbox{$\w=\ones,\, \g=\zeros$, 620 routes \label{fig:UniHbf-w1g0}}	[.35\linewidth]{\includegraphics[scale=0.1949]{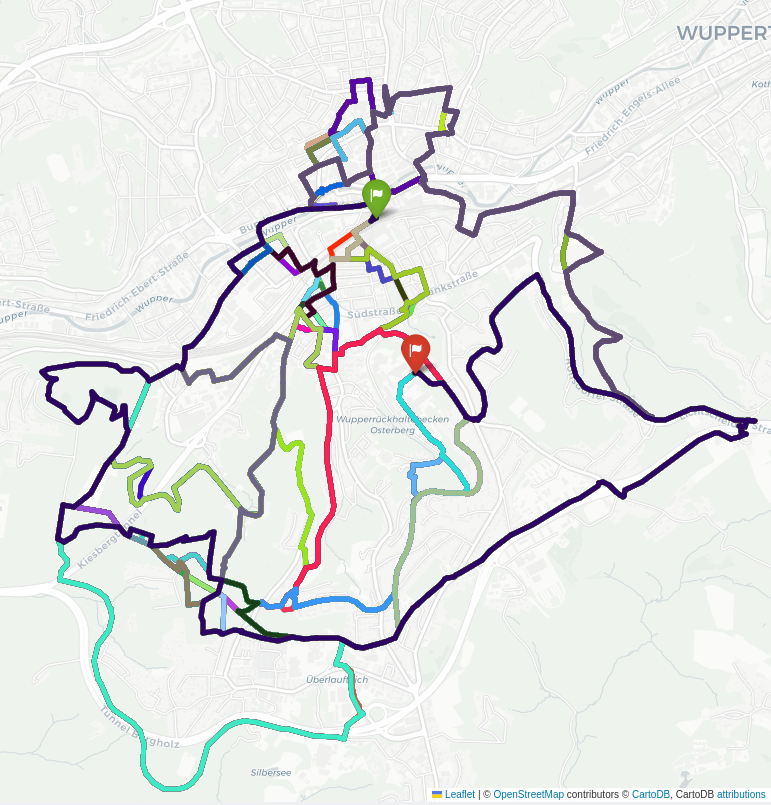}}
    \subcaptionbox{$\w=\ones,\, \g=0.2\cdot\ones$, 23 routes  }	[.28\linewidth]{\includegraphics[scale=0.3]{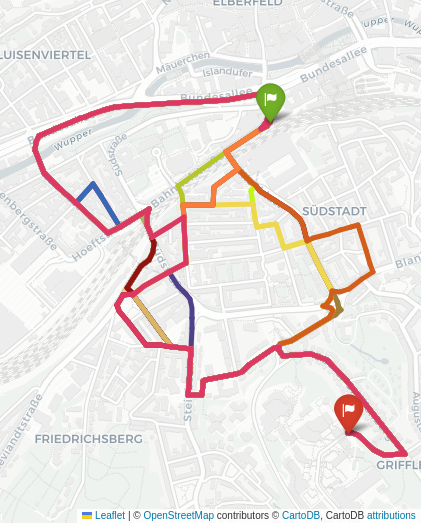}}
    \subcaptionbox{$\w=\ones,\, \g=0.4\cdot\ones$, 4 routes  }	[.28\linewidth]{\includegraphics[scale=0.3]{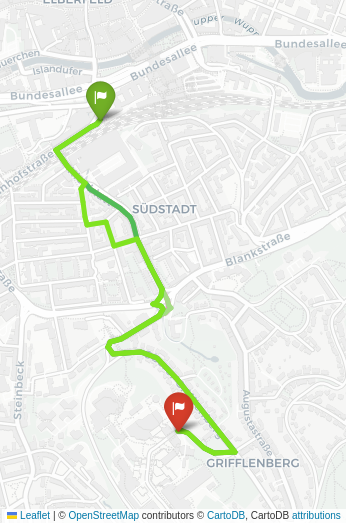}}\\
    \subcaptionbox{$\w=1.5\cdot\ones,\, \g=\zeros$, 193 routes  }	[.35\linewidth]{\includegraphics[scale=0.1949]{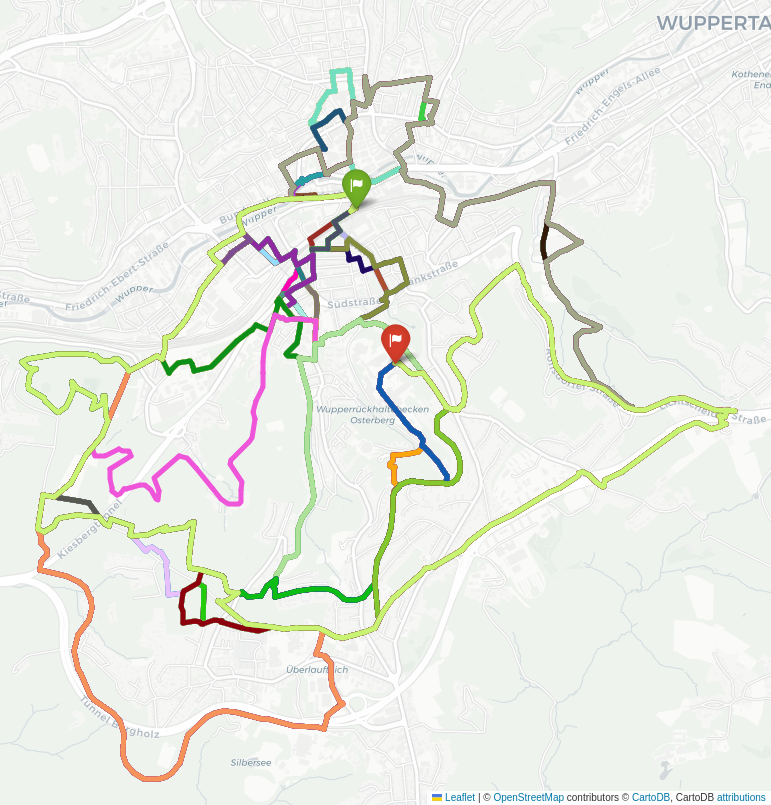}}
    \subcaptionbox{$\w=1.5\cdot\ones,\, \g=0.2\cdot\ones$, 21 routes  }	[.28\linewidth]{\includegraphics[scale=0.3]{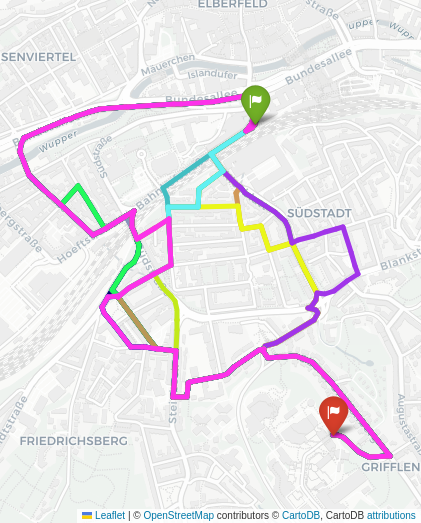}}
    \subcaptionbox{$\w=1.5\cdot\ones,\, \g=0.4\cdot\ones$, 1 route  }	[.28\linewidth]{\includegraphics[scale=0.3]{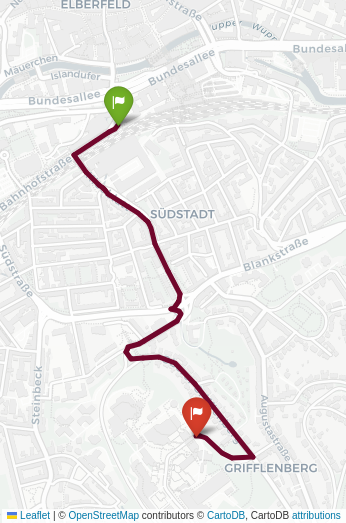}}\\
    \subcaptionbox{$\w=2\cdot\ones,\, \g=\zeros$, 131 routes  }	[.35\linewidth]{\includegraphics[scale=0.1949]{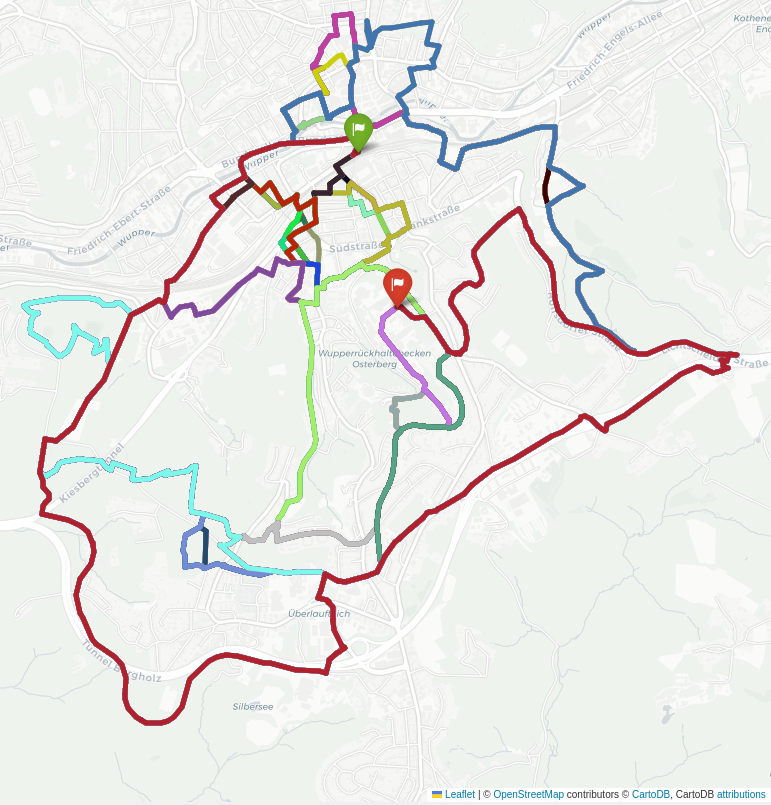}}
    \subcaptionbox{$\w=2\cdot\ones,\, \g=0.2\cdot\ones$, 16 routes  }	[.28\linewidth]{\includegraphics[scale=0.3]{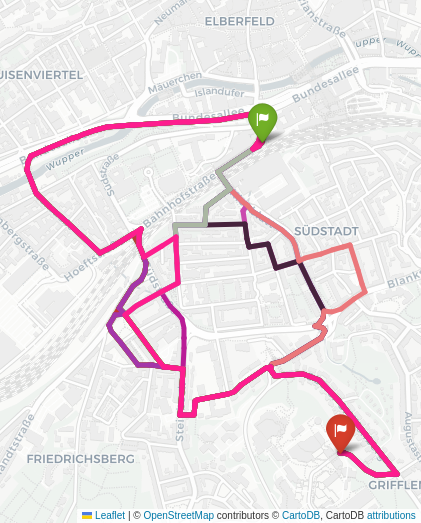}}
    \subcaptionbox{$\w=2\cdot\ones,\, \g=0.4\cdot\ones$, 1 route  }	[.28\linewidth]{\includegraphics[scale=0.3]{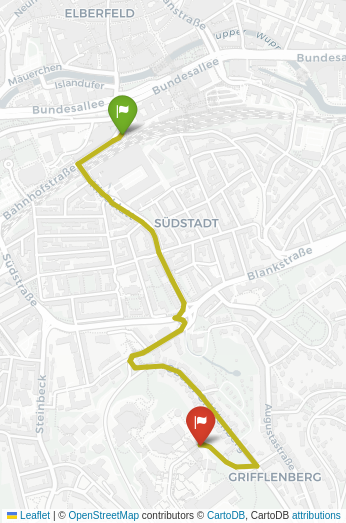}}\\
    \caption{Ordinal safest paths for different weights from the main campus of the University of Wuppertal to the main station of Wuppertal, based on \cite{OpenStreetMap}. Each path has a different color.}
	\label{fig:UniHbf}
\end{figure}

\begin{figure}
	\centering
    \subcaptionbox{$\w=\ones,\, \g=\zeros$, longest route \label{fig:UniHbf-longPath}}	[.35\linewidth]{\includegraphics[scale=0.1949]{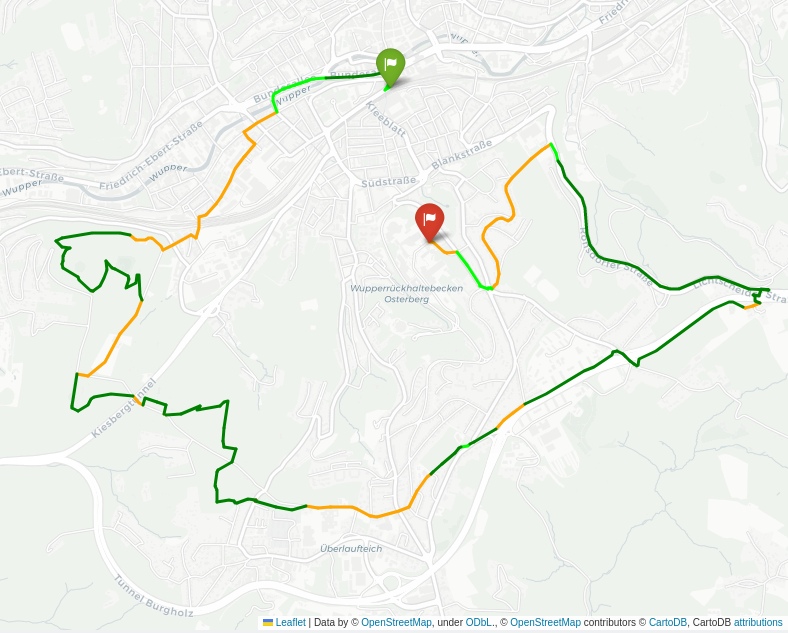}} 
    \subcaptionbox{$\w=\ones,\, \g=\zeros$, shortest route \label{fig:UniHbf-shortPath}}	[.28\linewidth]{\includegraphics[scale=0.15]{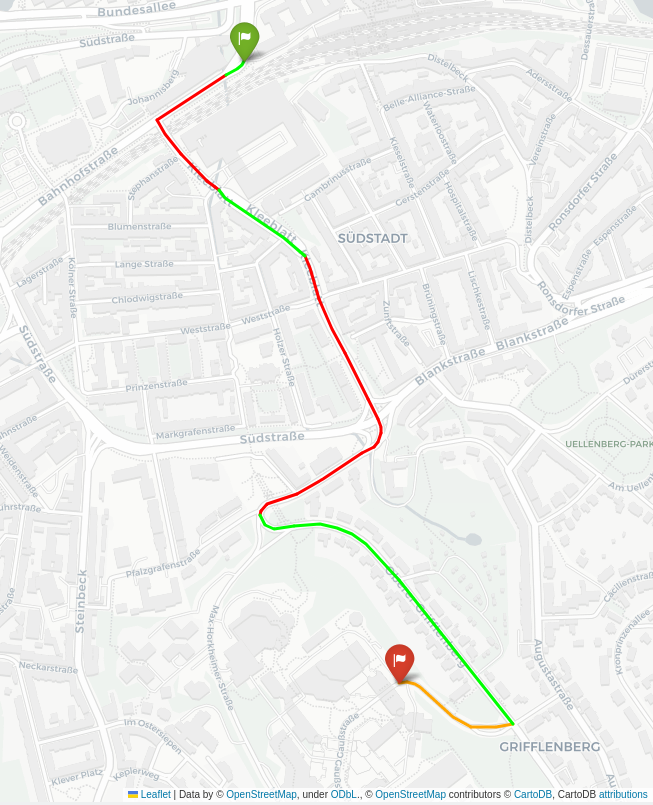}} 
    \subcaptionbox{$\w=2\cdot\ones,\, \g=0.4\cdot\ones$, only route \label{fig:UniHbf-onlyPath}}	[.28\linewidth]{\includegraphics[scale=0.15]{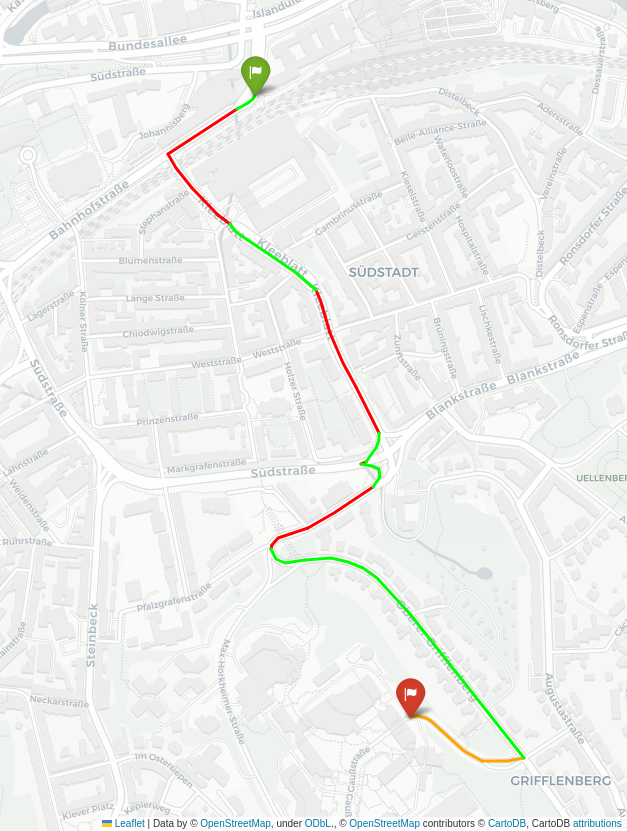}}
    \caption{Ordinally safest paths for different weights from the main campus of the University of Wuppertal to the main station of Wuppertal, based on \cite{OpenStreetMap}. The paths are colored regarding to their safety.}
	\label{fig:UniHbf-coloredPaths}
\end{figure}

We consider two different origin-destination pairs: first from the Campus Grifflenberg (main campus) of the University of Wuppertal to the main station of Wuppertal and second from the Campus Grifflenberg to Campus Haspel of the University of Wuppertal. 
In the first case we restrict the map to a square with the starting point as center and with a side length of $8$ km, i.\,e., $4$ km from the starting point to a side of the square, which seems reasonable as the straight line ($\ell_2$)-distance between the origin and the destination is approximately $1$ km. The number of ordinally non-dominated paths with $\omega=\ones$ and $\gamma=\zeros$ is $620$ on this map. All paths are depicted in Figure~\ref{fig:UniHbf-w1g0}, where each path has a different color. Note that these paths are partially overlapping and that the high number of paths results on the high number of possibilities to combine detours with each other. Moreover, one of the shortest paths and one of the longest ordinally non-dominated paths are visualized with their categories (dark green $=\eta_1$, light green $=\eta_2$, orange $=\eta_3$ and red $=\eta_4$) in Figures~\ref{fig:UniHbf-shortPath} and \ref{fig:UniHbf-longPath}, respectively. The shorter path takes many dangerous streets, while the longer path uses an enormous detour (which would not be considered in practice). 

To avoid very short but dangerous routes, we can increase the weights $\omega$. If we choose, for example, $\omega=1.5\cdot\ones$, i.\,e., we prefer one and a half meters of category $\eta_i$ over one meter of category $\eta_{i+1}$ for $i=1,\dots,3$, the number of efficient solutions reduces to $193$. If we further increase $\omega$, such that we prefer two meters of category $\eta_i$ over one meter of category $\eta_{i+1}$ for $i=1,2,3$, i.\,e., $\omega=2\cdot\ones$, the number of efficient solutions reduces to $131$. The sets of efficient solutions are visualized in Figure~\ref{fig:UniHbf} in the first column. We can see that the short routes and some specific subpaths are eliminated. Nonetheless, many very long routes are still efficient.

To eliminate those very long routes, we have to increase the weight of $\gamma$. Thus, we consider a weighted ordinal ordering cone $\tailC$ for $K=4$ with strictly positive weights. Then the corresponding matrix $A$ is given by
    \begin{equation*}
        A=\begin{pmatrix}
            1 & \w_1 & \w_1\,\w_2 & \w_1\,\w_2\,\w_3\\
            \g_3 & \g_3\,\w_1 & \g_3\,\w_1\,\w_2 &\w_1\,\w_2\\
            \g_2 & \g_2\w_1 & \w_1 & \w_1\,\w_3\\
            \g_2\,\g_3 & \g_2\,\g_3\,\w_1 & \g_3\,\w_1 & \w_1\\
            \g_1 & 1 & \w_2 & \w_2\,\w_3\\
            \g_1\,\g_3 & \g_3 & \g_3\,\w_2 & \w_2 \\
            \g_1\,\g_2 & \g_2 & 1 & \w_3 \\
            \g_1\,\g_2\,\g_3 & \g_2\,\g_3 & \g_3 & 1
        \end{pmatrix}.
    \end{equation*}
We consider for $\gamma$ the weights $\gamma=0.2\cdot\ones$ and $\gamma=0.4\cdot\ones$, respectively, i.\,e., we assume that one meter of category $\eta_{i+1}$ is preferred over $5$ meters or $2.5$ meters, respectively, in category $\eta_i$ for all $i=1,2,3$. The number of efficient solutions reduces significantly by increasing $\gamma$, see Figure~\ref{fig:UniHbf}, where $\omega$ is increased from the top to the bottom and $\gamma$ is increased from the left to the right. Obviously, the number of long routes is drastically reduced and for $\gamma=0.4\cdot\ones$ and $\omega\geqq 1.5\cdot\ones$ only one route is efficient, which is visualized with colors corresponding to the safety of the path segments in Figure~\ref{fig:UniHbf-onlyPath}. This route is quite similar to the shortest route for $\omega=\ones$ and $\gamma=\zeros$ in Figure~\ref{fig:UniHbf-shortPath}. It only differs in the crossing of a junction. However, the influence of the weights $\g,\w$ is clearly visible in this example.

\begin{figure}
	\centering
    \subcaptionbox{$\w=\ones,\, \g=\zeros$, $932$ routes }	[.3\linewidth]{\includegraphics[scale=0.22]{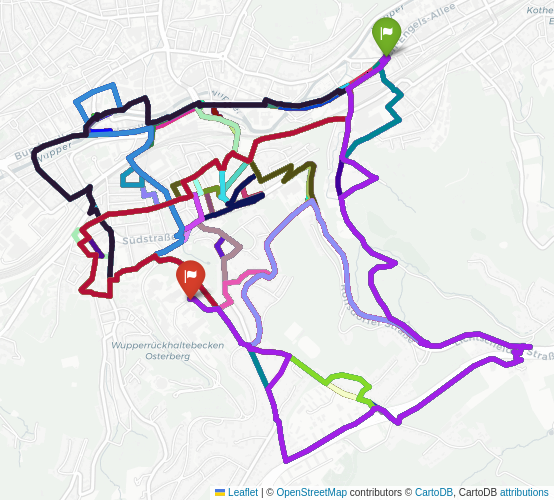}}
    \subcaptionbox{$\w=\ones,\, \g=0.2\cdot\ones$, $345$ routes  }	[.3\linewidth]{\includegraphics[scale=0.22]{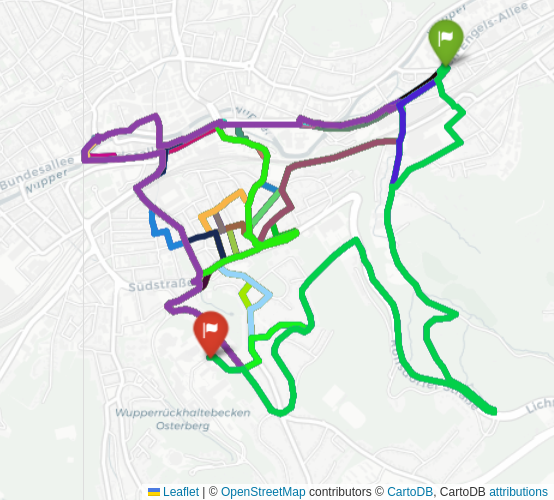}}
    \subcaptionbox{$\w=\ones,\, \g=0.4\cdot\ones$, $150$ routes  }	[.3\linewidth]{\includegraphics[scale=0.22]{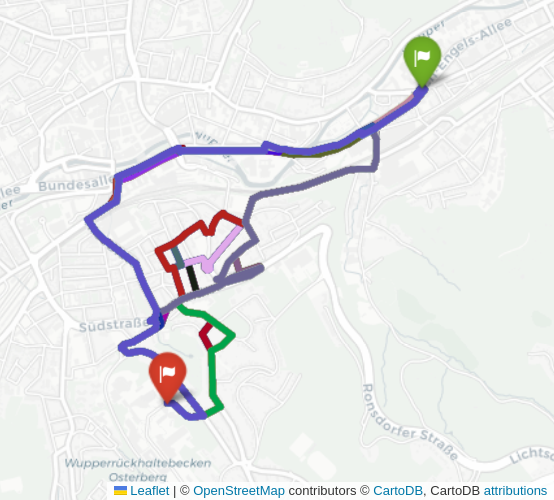}}\\
    \subcaptionbox{$\w=1.5\cdot\ones,\, \g=\zeros$, $141$ routes  }	[.3\linewidth]{\includegraphics[scale=0.22]{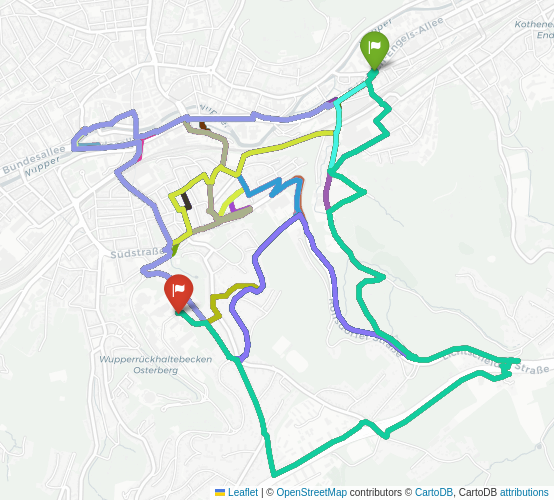}}
    \subcaptionbox{$\w=1.5\cdot\ones,\, \g=0.2\cdot\ones$, $65$ routes  }	[.3\linewidth]{\includegraphics[scale=0.22]{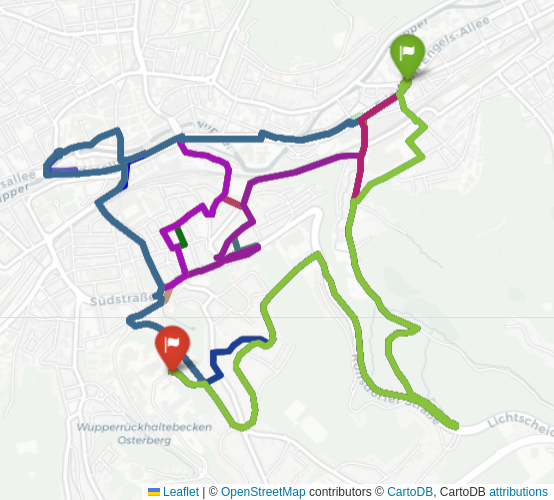}}
    \subcaptionbox{$\w=1.5\cdot\ones,\, \g=0.4\cdot\ones$, 24 routes  }	[.3\linewidth]{\includegraphics[scale=0.22]{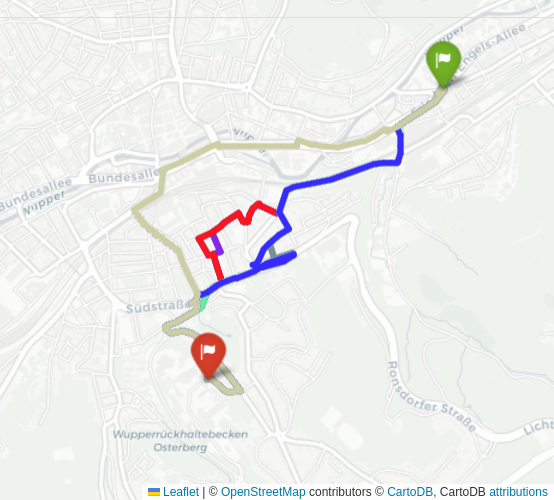}}\\
    \subcaptionbox{$\w=2\cdot\ones,\, \g=\zeros$, $27$ routes  }	[.3\linewidth]{\includegraphics[scale=0.22]{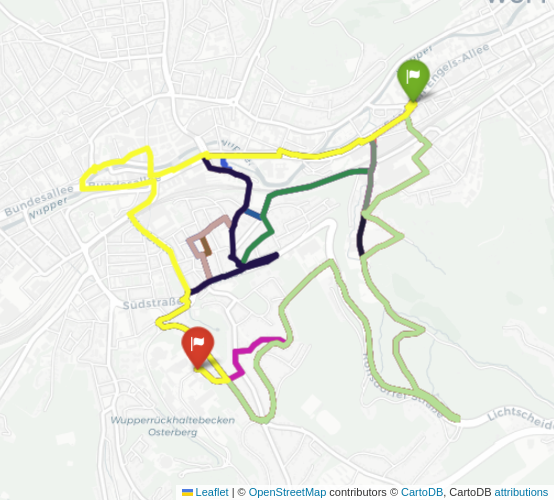}}
    \subcaptionbox{$\w=2\cdot\ones,\, \g=0.2\cdot\ones$, $14$ routes  }	[.3\linewidth]{\includegraphics[scale=0.22]{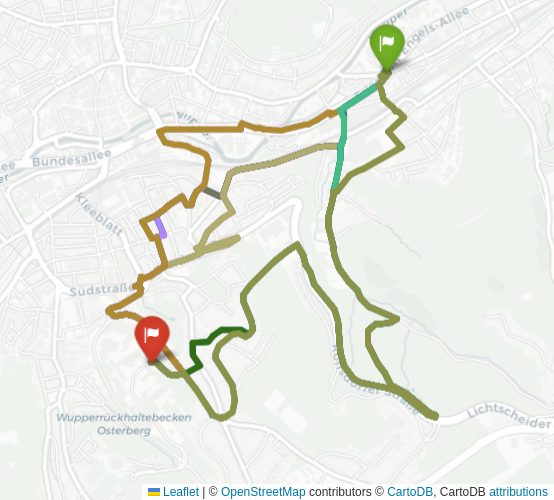}}
    \subcaptionbox{$\w=2\cdot\ones,\, \g=0.4\cdot\ones$, $5$ routes  }	[.3\linewidth]{\includegraphics[scale=0.22]{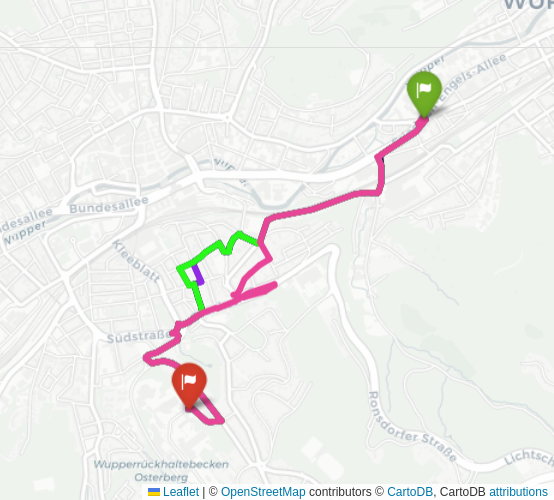}}\\
    \caption{Ordinal safest paths for different weights from Campus Grifflenberg to Campus Haspel of the University of Wuppertal, based on \cite{OpenStreetMap}.}
	\label{fig:GrifflenbergHaspel}
\end{figure}

\begin{figure}
	\centering
    \subcaptionbox{$\w=\ones,\, \g=\zeros$, longest route \label{fig:GrifflenbergHaspel-longPath}}	[.45\linewidth]{\includegraphics[scale=0.2]{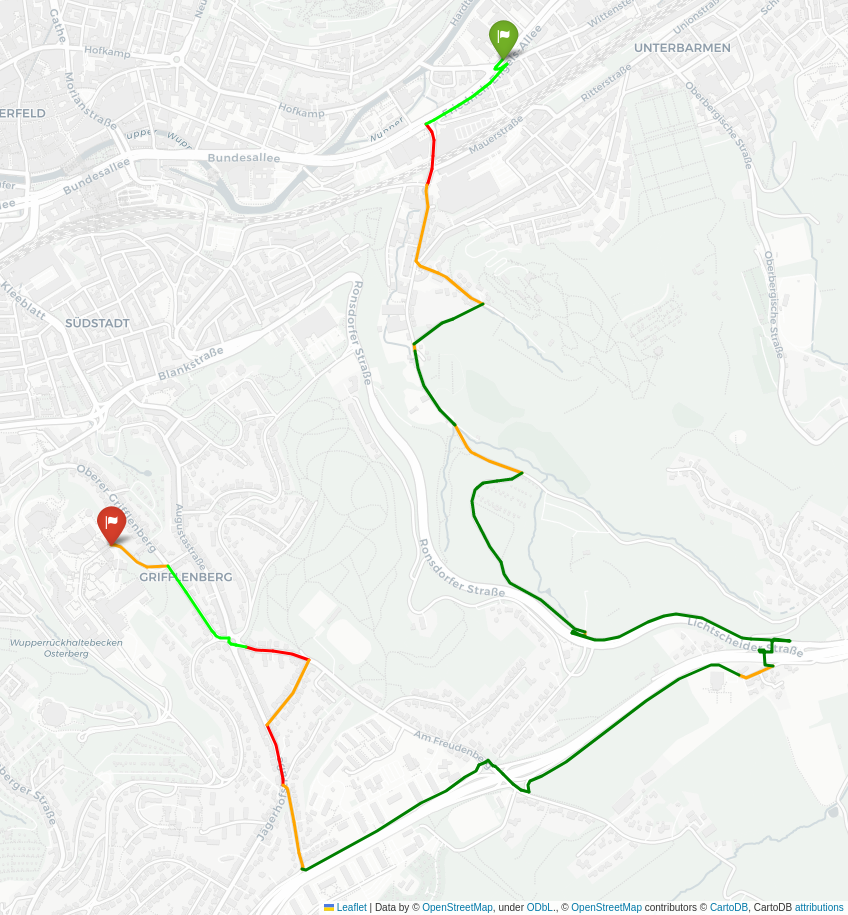}} 
    \subcaptionbox{$\w=\ones,\, \g=\zeros$, shortest route \label{fig:GrifflenbergHaspel-shortPath}}	[.45\linewidth]{\includegraphics[scale=0.29]{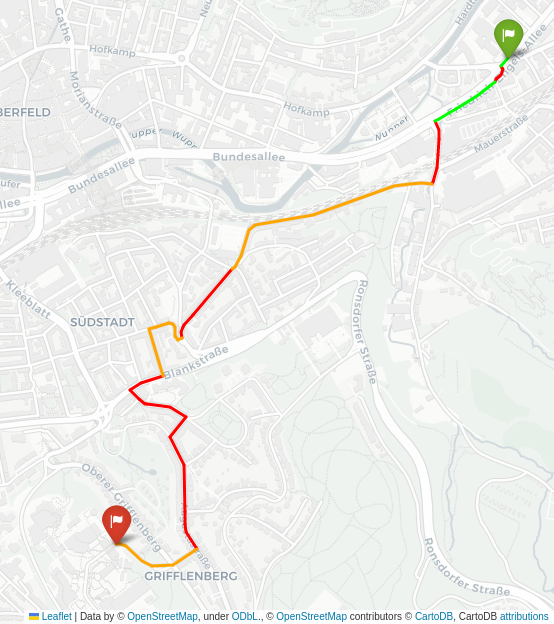}} 
    \caption{Ordinal safest paths for different weights from Campus Grifflenberg to Campus Haspel of the University of Wuppertal, based on \cite{OpenStreetMap}. The paths are colored regarding to their safety.}
	\label{fig:GrifflenbergHaspel-w1g0-coloredPaths}
\end{figure}

We consider a second origin-destination pair from the Campus Grifflenberg to Campus Haspel of the University of Wuppertal to show that the reduction of the alternative routes by increasing $\omega$ and $\gamma$ is highly dependent on the instance. In this case, we restrict the map to a square with the starting point as center and with a side length of $6$ km. Again, we considered different values for $\omega\in\{\ones,\, 1.5\cdot\ones,\, 2\cdot\ones\}$ and $\gamma\in\{\zeros,\,0.2\cdot\ones,\,0.4\cdot\ones\}$ and the resulting efficient paths are visualized in Figure~\ref{fig:GrifflenbergHaspel}. Here, the reduction of the number of routes is less significant when increasing $\gamma$. However, for $\omega=2\cdot\ones$ and $\gamma=0.4\cdot\ones$ there exist only $5$ efficient paths. Moreover, there are only three paths that use different streets, while the other only differ on a junction (visualized in Figure~\ref{fig:GrifflenbergHaspel-w2-g0,4-coloredpaths}). These paths are compromise solutions and do not correspond to the longest or shortest path for $\omega=\ones$ and $\gamma=\zeros$, which are depicted in Figure~\ref{fig:GrifflenbergHaspel-w1g0-coloredPaths}.

\begin{figure}
	\centering
    \subcaptionbox{Route 1 }	[.3\linewidth]{\includegraphics[scale=0.22]{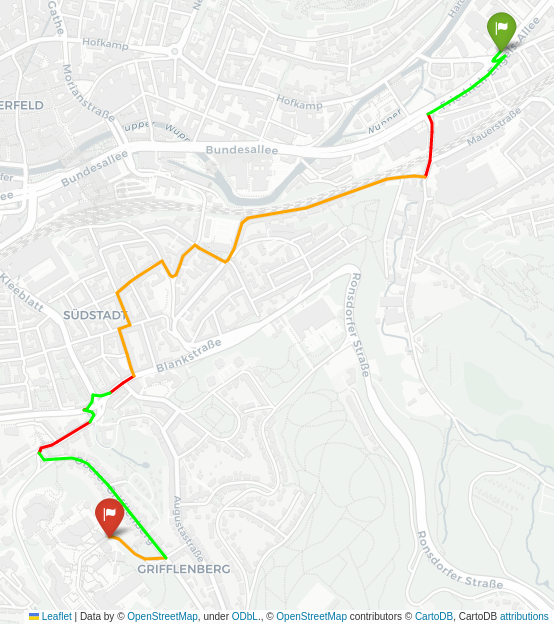}}
\subcaptionbox{Route 2 }	[.3\linewidth]{\includegraphics[scale=0.22]{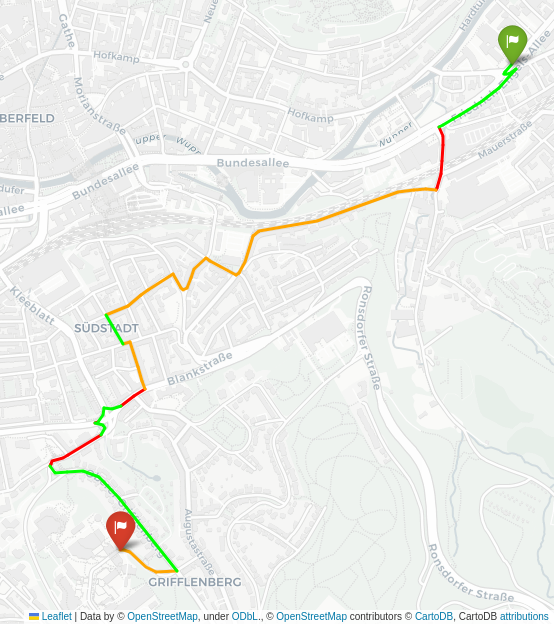}} \subcaptionbox{Route 3  }	[.3\linewidth]{\includegraphics[scale=0.22]{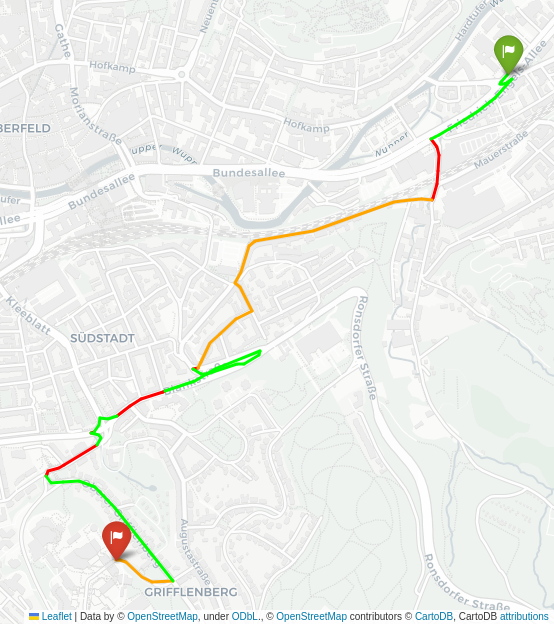}}
    \caption{Ordinal safest paths for $\w=2\cdot\ones,\, \g=0.4\cdot\ones$ from Campus Grifflenberg to Campus Haspel of the University of Wuppertal, based on \cite{OpenStreetMap}. The paths are colored regarding to their safety.}
	\label{fig:GrifflenbergHaspel-w2-g0,4-coloredpaths}
\end{figure}

\section{Conclusion}\label{sec:conclusion}

In this article, we investigate ordinal combinatorial optimization problems and introduce marginal weights for the categories. The weights are used to model additional preferences between consecutive categories, so that only relevant solutions are computed. We first introduce the concept of standard ordinal optimality with $K$ ordered categories and extend this to weighted ordinal optimality.

We show that weighted ordinal dominance can be modeled by a closed, polyhedral ordering cone with up to $2(K-1)$ extreme rays and $2^{K-1}$ facets. We present the corresponding matrices that describe this cone by its extreme rays and by its facets. The latter can be used to transform the problem linearly to an associated multi-objective optimization problem with as many objective functions as the weighted ordinal ordering cone has facets.

We analyze some specific cases, where the number of objective functions does not change through the transformation, i.\,e., the number of facets of the weighted ordinal ordering cone equals the number of categories. Furthermore, we interrelate the weighted ordinal ordering cone with other cones which are often used in multi-objective optimization like, e.\,g., the Pareto cone, the lexicographic ordering cone, or the standard ordinal cone.

The influence of the weights on the efficient solutions is visualized by computing safe bike routes in Wuppertal between the main campus of the University  and the main station or the campus Haspel using different weights, respectively. It is shown that without weights it may happen that some very short but dangerous routes as well as extremely long routes are computed. The weights may reduce the cardinality of the set of efficient solutions drastically and reduce it to relevant solutions only.

The case study also shows that considering only weights on the streets is in general not sufficient, as crossing the street (by foot) or turning left or right (by bike) is related to risk that is not reflected by the edge. Thus, in future work, one could analyze how the crossing points can be modeled more realistically.

A further extension of our approach could be to introduce weights not only between consecutive categories, but between all pairs of categories. This generalizes the concept of weighted ordinal optimization and is step towards general polyhedral dominance cones. However, in this situation it is very difficult to chose weights that represent a well-defined, non-conflicting preference structure.

\subsection*{Acknowledgements}\noindent
This work is partially supported by the  project SAFER -- Safe, Algorithm-based Footpath Routing, which is funded by the EU program EFRE/JTF NeueWege.IN.NRW (EFRE 20800941).

\end{document}